\documentclass[final,onefignum,onetabnum]{article}

\usepackage{array}

\usepackage[caption=false]{subfig}
\usepackage{pgfplots}

\newtheorem{thm}{Theorem}[section]

\newtheorem{rem}[thm]{Remark}

\usepackage{algorithmic}

\usepackage{graphicx,epstopdf}

\usepackage{subfig}
\usepackage{amsmath}
\usepackage{amssymb}
\usepackage{graphicx}
\usepackage{dcolumn}
\usepackage{mathtools}
\usepackage{amsmath,amsfonts}
\usepackage{bm}
\usepackage{amsmath}
\usepackage{amssymb}
\usepackage{color}
\usepackage{float}
\usepackage{setspace}
\usepackage{tabularx}

\allowdisplaybreaks

\newcommand{\ds}{\displaystyle}
\newcommand{\E}{ {\mathbb{E}} }

\newcommand{\be}{\begin{equation}}
\newcommand{\ee}{\end{equation}}

\def\cF{{\cal F}} \def\dbF{\mathbb{F}} \def\dbP{\mathbb{P}} \def\ges{\geqslant} \def\les{\leqslant} \def\qq{\qquad}
\def\q{\quad} \def\dbR{\mathbb{R}} \def\cU{{\cal U}}  \def\hb{\hbox} \def\rf{\eqref}
 \def\cD{{\cal D}}

\def\lan{\langle} \def\ran{\rangle}

              \def\Om{\Omega}   
                    
\def\z{\zeta}               \def\si{\sigma}
\def\e{\varepsilon}             
              
\def\ba{\begin{array}}                \def\ea{\end{array}}
\def\bel{\begin{equation}\label}      \def\ee{\end{equation}}

\colorlet{texcscolor}{blue!50!black}
\colorlet{texemcolor}{red!70!black}
\colorlet{texpreamble}{red!70!black}
\colorlet{codebackground}{black!25!white!25}

\date{}

\title{An efficient numerical algorithm for solving data driven feedback control problems}

\author{
Richard Archibald\thanks{ Computational Science and Mathematics Division, Oak Ridge National Laboratory, Oak Ridge, Tennessee.}
\and Feng Bao\thanks{ Department of Mathematics, Florida State University, Tallahassee, Florida, \ ({\tt bao@math.fsu.edu}).}
\and Jiongmin Yong\thanks{ Department of Mathematics, University of Central Florida, Orlando, Florida,}
\and Tao Zhou\thanks{ Academy of Mathematics and Systems Sciences, Chinese Academy of Sciences, Beijing, China, }
      }

\begin{document}
\maketitle

\begin{abstract}
  The goal of this paper is to solve a class of stochastic optimal control problems numerically, in which the state process is governed by an It\^o type stochastic differential equation with control process entering both in the drift and the diffusion, and is observed partially. The optimal control of feedback form is determined based on the available observational data. We call this type of control problems the data driven feedback control. The computational framework that we introduce to solve such type of problems aims to find the best estimate for the optimal control as a conditional expectation given the observational information. To make our method feasible in providing timely feedback to the controlled system from data, we develop an efficient stochastic optimization algorithm to implement our computational framework.
\end{abstract}


\textbf{AMS:} 93E11, 60G35, 65K10	

\section{Introduction}

Stochastic optimal control is an important research subject that attracts scientists and engineers in various fields from theoretical scientific research to practical industrial production. The {\it control process} (also called {\it control policy}), which controls a stochastic dynamical system whose solution called {\it state process}, is designed to meet some optimality conditions. For the classic stochastic optimal control problem with full observation of the state, both theoretical results and numerical methods are extensively studied. However, in practice the full observation of the state process is often not available. Instead, we have detectors/observation facilities to collect partial observational data, which provide indirect information about the state process. The theoretical formulation for partially observable stochastic optimal control is derived analytically \cite{Fleming-Pardoux-1982, Haussmann-1987, Li-Tang-1995, Tang-1998, Charalambous-98, Wang-Wu-Xiong-2018}, and the corresponding optimal control is a stochastic process adapted to the observational information. Since the controlled state needs to be inferred from observations, the procedure of finding the optimal control requires data analysis for observational data, and the control actions are driven by information contained in data. To highlight the influence of data in designing control policies, we call the stochastic optimal control with partially observed controlled processes the {\it data driven feedback control}.

The key to solve a data driven feedback control problem is the effective combination of data analysis techniques and stochastic optimal control methods.
One of the most widely used methods to combine data with control is the ``separation principle'' \cite{Anders_Separation, Wonham_Separation}. The main theme of separation principle is that the data driven feedback control problem can be separated into a state estimation problem, and a standard stochastic optimal control problem, in which the state of controlled process is provided by the state estimation. Despite of its significant success in practical applications, the separation principle only works for time-invariant linear controlled systems with linear observations. Although several efforts have been made to extend the applicability of the separation principle by linearizing the nonlinear models \cite{Charalambous-98, Linearization_EKF, Local_Separation},  computational methods for general \textit{nonlinear} problems are still needed.


The goal of this paper is to introduce an effective computational framework that incorporates observational information into stochastic optimal control problems and we shall also design an efficient algorithm to implement our computational framework through stochastic optimization.
Since the data driven feedback control is adapted to observational information, we use the conditional expectation for the control process as its ``best'' approximation given observations. In this way, our computational framework is composed by two components: (I) numerical schemes to obtain conditional distribution of the controlled process given observations; (II) computational methods for stochastic optimal control problems. Component (II) provides the overarching algorithmic structure to find the optimal control, and the conditional distribution obtained in component (I) is used to calculate conditional expectations for the optimal control solver.
In this paper, we apply optimal filtering methods, which could find the conditional distribution of a dynamical system based on indirect observations \cite{BSDE_filter}, to accomplish component (I) in our computational framework. When solving nonlinear optimal filtering problems, Zakai filter and particle filter are two well-known approaches. The methodology of the Zakai filter is to solve for the conditional probability density function (pdf) for the target dynamical system through a parabolic type stochastic partial differential equation (called Zakai equation) \cite{zakai}; and the particle filter, also known as a sequential Monte Carlo method, describes the desired conditional pdf by using empirical distribution of a set of random samples (particles) \cite{particle-filter}. Although the Zakai filter could provide more accurate approximation for conditional distributions theoretcally, the particle filter is more popular in solving practical problems due to the high efficiency of Monte Carlo method in approximating high dimensional distributions \cite{MTAC2012}. Therefore, in our approach we use the particle filter to approximate the conditional distribution of the controlled process.
On the other hand, computational methods for solving the stochastic optimal control problem in component (II) are developed under two general frameworks -- dynamic programming principle and maximum principle \cite{Yong-Zhou-199, Feng_2013, DG_HJB_2016, Gong_2017}. In this work, we adopt the maximum principle framework to solve the optimal control problem due to its multiple advantages over dynamic programming principle. For example, it has less restrictions on dimension of the problem and it could solve problems with state constraints -- especially with finite dimensional terminal state constraints, and it also allows to have random coefficients in the controlled state equation and/or in the performance cost functional. In maximum principle, we have a stochastic Hamiltonian system that consists of a system of forward backward stochastic differential equations (FBSDEs), which meet certain optimality condition with respect to the optimal control \cite{Peng_control, Touzi_control}. Therefore, solving the stochastic optimal control problem through maximum principle involves solving the FBSDEs system (usually carried out by numerical methods) and an optimization procedure, which is typically achieved by gradient descent iterations.
There are several successful numerical methods for solving FBSDEs, which can be categorized by two types of approaches: the first  type of approach solves an FBSDEs system through numerical schemes for its equivalent parabolic partial differential equation \cite{Four_step, Milstein_BSDE}; the second approach solves FBSDEs as a system of stochastic differential equations and solutions of FBSDEs are approximated by their conditional expectations \cite{Bally_BSDE, Touzi_BSDE,  Exarchos_Control, Tsiotras_control, Zhao_BSDE}. While both approaches effectively solve FBSDEs, implementing numerical algorithms to solve FBSDEs is still a challenging task -- even with some advanced high efficiency methods \cite{ML_BSDE},  and the computational cost for solving FBSDEs increases dramatically when the dimension of the problem becomes higher.  In our approach, we choose the second approach to solve FBSDEs due to the efficiency and flexibility of stochastic computing methods.

From the above discussion, we can see that the implementation of our data driven feedback control framework includes particle filtering for conditional pdf of the controlled precess and numerical optimization for control process which involves numerical solutions for FBSDEs. Since the purpose of data driven feedback control is to give timely feedback to the controlled system based on observational data, efficiency of a numerical method is essential. To provide a highly efficient algorithm to implement our computational framework, we introduce a stochastic optimization algorithm that combines the particle filter with stochastic gradient descent in the optimization procedure for the data driven optimal control.  The motivation of applying stochastic optimization is the fact that our computational framework consists multiple conditional expectations and stochastic gradient descent is a very efficient method to treat expectations in the gradient descent type optimization. In this connection, we consider all the random samples, including particles obtained from the particle filter, to be the ``data'' which we use to calculate conditional expectations. Then, in each gradient descent iteration step, we pick one random sample in the sample space, as well as one particle in the particle cloud, and represent conditional expectations by the realization of stochastic process generated by our picked sample/particle. In this way, the gradient decent becomes its stochastic counterpart and we transfer the computational cost of fully calculating conditional expectations (by averaging all samples/particles) to more iteration steps with one-sample simulation of stochastic process in each step. The justification of stochastic gradient descent is well studied theoretically \cite{SGD-E-2019, Convergence-SGD} and the efficiency is verified in many application scenarios.
The only potential issue for our application of stochastic optimization is that the single-realization representation of conditional expectations could not solve FBSDEs accurately.
Indeed, the usage of only one realization of random sample is not sufficient to provide full characterization for solutions of FBSDEs. However, the purpose of solving the FBSDEs system in maximum principle is to formulate a stochastic Hamiltonian system, which is used to describe the gradient process with respect to the control process. Therefore solving FBSDEs is not the goal of our computational framework. Actually, the conditional expectations that we use to approximate solutions of FBSDEs are eventually used to describe the gradient process, and similar justifications for stochastic gradient descent would also apply to FBSDEs in our approach.
Moreover, by using the single-realization representation of expectation, we can avoid unnecessary calculations for obtaining fully calculated solutions of FBSDEs, which is the primary computational cost of maximum principle approach for stochastic optimal control problems.

The rest of this paper is organized as following. In Section \ref{Theory}, we briefly discuss the theoretical background for the data driven feedback control problem and introduce an optimization framework for finding the optimal control, which will be used to design our numerical algorithms. The computational framework with our efficient stochastic optimization implementation will be introduced in Section \ref{Numerical approach}. In Section \ref{Numerics}, we demonstrate the effectiveness and efficiency of our numerical algorithms by solving both classic benchmark problems and practical feedback control problems. Some concluding remarks and the plan of our future work are given in Section \ref{Conclusion}

\section{Data driven feedback control problem}\label{Theory}

\subsection{Problem statement}\label{DDFC}

Let $(\Om,\cF,\dbF^W,\dbP)$ be a complete filtered probability space on which a $d$-dimensional standard Brownian motion $W:=\{W_s\}_{s \ges 0 }$ is defined with $\dbF^W=\{\cF^W_s\}_{s\ges0}$ being its natural filtration augmented by all the $\dbP$-null sets in $\cF$. Consider the following stochastic differential equation (SDE) over a deterministic time interval $[0,T]$
\bel{state[0,T]}
dX_t=b(t,X_t,u_t )dt+\si(t,X_t,u_t)dW_t,\qq t\in[0,T],\qq
X_0=\xi,\ee
where $X$ is the (controlled) state process valued in $\dbR^d$, $u$ is the control process valued in some set $U\subseteq\dbR^m$, $b:[0,T]\times\dbR^d\times U\to\dbR^d$ and $\si:[0,T]\times\dbR^d\times U\to\dbR^{d\times d}$ are suitable maps, and $\xi$ is a random variable independent of $W$ following a distribution $p_0$. We call the terms $b$ and $\si$ the {\it drift} and the {\it diffusion}, respectively.
Let
$$\cU[0,T]=\big\{u:[0,T]\times\Om\to U\bigm|u \hb{ is $\dbF^W$-progressively measurable }\big\}.$$
Under some mild conditions, for every square integrable initial variable $\xi$ and control $u\in\cU[0,T]$, equation \eqref{state[0,T]} admits a unique solution $X$. To measure the performance of the control $u$, we introduce the following cost functional
\bel{cost[t,T]}J(u)=\E\left[ \int_0^Tf(t,X_t,u_t)dt+h(X_T)\right].\ee
Then classical optimal control problem can be stated as

\vspace{0.5em}
\bf Problem (C). \rm For any given initial condition $X_0 = \xi$, find a $\bar u \in\cU[0,T]$ such that
\bel{J=inf}
J(\bar u)=\inf_{u\in\cU[0,T]}J(u).
\ee
Any $\bar u\in\cU[0,T]$ satisfying the above is called an {\it optimal control}. 
Under proper conditions, if the optimal control $\bar u$ exists, it should be a function of the corresponding state process $\bar X$, in some sense. However, in reality, we may not be able to observe the true state $X$. Instead, we introduce the following equation
\bel{observe[0,T]}dM_t=g(X_t)dt+dB_t, \qq M_0 = 0, \ee
where $g: \dbR^d\to\dbR^\ell$ and $B$ is an $\ell$-dimensional standard Brownian motion independent of $W$. We call the above an observation equation and $M$ is called the {\it observation process}, which provides partial (noisy) observations for $X$. We may put \eqref{state[0,T]} and \eqref{observe[0,T]} together to get the following augmented system on time interval $[0,T]$
%
%
%
\bel{XM}
\ds d\begin{pmatrix}X_t\\ M_t\end{pmatrix} =\begin{pmatrix}b(t,X_t,u_t )\\ g(X_t)\end{pmatrix}dt  + \begin{pmatrix}\sigma(t,X_t,u_t)&0\\0&I\end{pmatrix}
d\begin{pmatrix}W_t\\ B_t\end{pmatrix}, \qq
\begin{pmatrix}X_0 = \xi \\ M_0 = 0 \end{pmatrix}.
\ee
Let $\dbF^B=\{\cF_t^B\}_{t\ges0}$ be the filtration of $B$ augmented by all the $\dbP$-null sets in $\cF$, and $\dbF^{W,B}\equiv\{\cF^{W,B}_t\}_{t\ges0}$ be the filtration generated by $W$ and $B$ (augmented by $\dbP$-null sets in $\cF$).
Under mild conditions, for any square integrable random variable $\xi$ independent of $W$ and $B$, and any $\dbF^{W,B}$-progressively measurable process $u$ (valued in $U$), \rf{XM} admits a unique solution $(X,M)$ which is $\dbF^{W,B}$-adapted. Next, we let $\dbF^M=\{\cF^M_t\}_{t\ges0}$ be the filtration generated by $M$ (augmented by all the $\dbP$-null sets in $\cF$). Clearly, $\dbF^M\subset\dbF^{W,B}$, and $\dbF^M\ne\dbF^W $, $\dbF^M\ne\dbF^B$, in general.
%
%
We introduce the set of admissible controls as
$$\cU_{ad}[0,T]=\big\{u:[0,T]\times\Om\to U\bigm|u\hb{ is $\dbF^M$-progressively measurable}\big\}.$$
To emphasis the  $\dbF^M$-adaptedness, for a control in the admissible control set $\cU_{ad}[0,T]$, we denote it by $u^M$. Replacing $u$ in \eqref{XM} by $u^M$ , we still have the unique solution $(X,M)$, which is $\dbF^{W,B}$-adapted. The cost functional should now be modified as
\bel{cost*[0,T]}J^*(u^M)=\E\Big[\int_0^Tf(t,X_t,u^M_t)dr+h(X_T)\Big].\ee
Then we pose the following optimal control problem.

\vspace{0.5em}
\bf Problem (C$^*$). \rm For any given square integrable random variable $\xi$, find a $u^*\in\cU_{ad}[0,T]$ such that
\bel{J^*=inf}J^*(u^{*})=\inf_{u^M\in\cU_{ad}[0,T]}J^*(u^M)\equiv V^*(t,\xi).\ee
It is worthy of mentioning that under the induced probability measure $P^M$ with Girsanov transformation $dP^M = \Theta_t^T dP$, where $\Theta_t^T :=\exp \big[- \int_t^T g(X_s) dB_s - \int_t^T \frac{1}{2} |g(X_s)|^2 ds \big]$, the observation process $M$ is a standard Brownian motion.
Discussions on the above indirect observations based stochastic optimal control can be traced back to the work of Fleming in 1960s \cite{Fleming-1968}, followed by Kwakernaak \cite{Kwakernaak-1981}, Fleming--Pardoux \cite{Fleming-Pardoux-1982}, Haussmann \cite{Haussmann-1982, Haussmann-1987}, Bensoussan \cite{Bensoussan-1983, Bensoussan-1992},  Baras--Elliott--Kohlmann \cite{Baras-Elliott-Kohlmann-1989}, Li--Tang \cite{Li-Tang-1995}, Tang \cite{Tang-1998}, and Wang--Wu--Xiong \cite{Wang-Wu-Xiong-2018}, and so on. Since the optimal control $u^*$ defined in \eqref{J^*=inf} implements feedback control actions based on the observational data, we name the optimal control problem (C$^*$) with its corresponding stochastic system \eqref{XM} \textit{``data driven feedback control problem''}.

\subsection{Optimization for optimal control}

In this work, we solve the data driven feedback control problem through an optimization procedure, which is derived by stochastic maximum principle. In the case that the optimal control $u^*$ is in the interior of $\mathcal{U}_{ad}$, one can deduce by using the G\^ateaux derivative of $u^*$ and maximum principle that the gradient process of the cost functional $J^*$ with respect to the control process on time interval $t \in [0, T]$ has the following form (see Appendix for details)
\begin{equation}\label{J'}
(J^*)'_u(u^*_t)=\E\big[b_u(t, X^*_t, u^*_t)^\top Y_t+\si_u(t, X^*_t, u^*_t)^\top Z_t+f_u(t, X^*_t, u^*_t)^\top\bigm|\cF^M_t\big],
\end{equation}
where we use subscripts to denote partial derivatives of functions throughout of this paper, and stochastic processes $Y$ and $Z$ are solutions of the following forward backward stochastic differential equations (FBSDEs) system
\bel{FBSDE}\left\{\ba{ll}
dX^*_t=b(t,X^*_t,u^*_t)dt+\si(t,X^*_t,u^*_t)dW_t,  \qq X_0 = \xi  & \text{(SDE)}\\
dM^*_t=g(X^*_t)dt+dB_t, \qq \qq \qq \qq \q M_0 = 0 & \text{(SDE)}\\
dY_t=\Big(-b_x(t,X^*_t,u^*_t)^\top Y_t - \si_x(t,X^*_t,u^*_t)^\top Z_t-f_x(t,X^*_t,u^*_t)^\top\Big)dt \\
\qq\qq\qq\qq  +Z_tdW_t+\z_tdB_t,  \qq Y_T=h_x(X^*_T)^\top,  &\text{(BSDE)}
\ea\right.
\ee
where $Z$ is the martingale representation of $Y$ with respect to $W$ and $\z$ is the martingale representation of $Y$ with respect to $B$. Due to the backward propagation direction, we call the third equation in \eqref{FBSDE} the backward stochastic differential equation (BSDE), which is derived as the \textit{adjoint equation} of the controlled state equation.

The main theme of our approach is to solve for the data driven feedback optimal control $u^*$ through gradient descent type optimization and the gradient process $(J^{*})'_{u}$  is defined in \eqref{J'}.
Specifically, for a pre-chosen $\dbF^M$-adapted process $u^{0,M}$ (as our initial guess), we introduce the following gradient descent iteration to find the optimal control $u^*_t$ at any time instant $t \in [0, T]$
\begin{equation}\label{GD}
u_t^{l+1, M} = u_t^{l,M} - \rho (J^{*})'_{u}(u_t^{l, M}), \quad l = 0, 1, \cdots,
\end{equation}
where $\rho$ is the step-size for the gradient, and we use $(J^{*})'_{u}(u)$ to denote the gradient process corresponding to a control process $u$. The primary challenge in applying the gradient descent optimization \eqref{GD} to solve for the data driven optimal control $u^*$ is to obtain the gradient process $\{(J^{*})'_{u}(u^{l, M}_t)\}_{0\les t\les T}$, which is described by solutions $Y_t$ and $Z_t$ of the FBSDEs system \eqref{FBSDE} corresponding to the $\cF^M_t$ adapted estimated optimal control $u^{l, M}_t$.
However, the observational information $\cF^M_t$ is gradually increased as we collect more and more data over time. At a certain time instant $t$, the information $\cF^M_s$ for $s \in [t, T]$ is not available to solve for $\{(Y_s, Z_s)\}_{t \les s \les T}$, and the estimated feedback control $\{u^{l, M}_s\}_{t \les s \les T}$ should not be driven by information $\{\cF^M_s\}_{t \les s \les T}$. Therefore, we target on finding the optimal control $u^{*}_t$ at the current time $t$ with accessible information $\cF^M_t$. Since the evaluation for $(J^*)'_u(u_t^{l,M})$ requires trajectories $\{(Y_s, Z_s)\}_{t \les s \les T}$ as $Y_t$ and $Z_t$ are solved backwards from $T$ to $t$, we take conditional expectation $\E[\,\cdot\,| \cF^M_t]$ to the gradient process $\{(J^*)'_u(u_s^{l,M})\}_{t \les s \les T}$, i.e.
\begin{equation}\label{J'_t}
\begin{aligned}
& \E[(J^*)'_u(u_s^{l,M}) | \cF^M_t] =  \E\Big[b_u(s, X_s, u_s^{l, M})^\top Y_s  \\
&\qq \qq \qq \qq +\si_u(s, X_s, u_s^{l, M})^\top Z_s+f_u(s, X_s, u_s^{l, M})^\top\bigm|\cF^M_t\Big], \quad  s \in[t, T],
\end{aligned}
\end{equation}
where $X_s$, $Y_s$ and $Z_s$ are corresponding to the estimated control $u_s^{l, M}$.
Similarly, for the gradient descent iteration \eqref{GD} on the time interval $[t, T]$, by taking  conditional expectation $\E[\,\cdot\,| \cF^M_t]$, we obtain
\begin{equation}\label{E-GD_s}
\E[u_s^{l+1, M} | \cF^M_t] = \E[ u_s^{l,M} | \cF^M_t] - \rho \E\big[ (J^{*})'_{u}(u_s^{l,M}) | \cF^M_t\big], \quad  l = 0, 1, \cdots,   \quad s \in [t, T].
\end{equation}
When $s=t$, the conditional expectation $\E[u_s^{l+1, M} | \cF^M_t]$ on the left hand side of \eqref{E-GD_s} gives us $u_t^{l+1, M} = \E[u_t^{l+1, M} | \cF^M_t]$ due to the $\cF^M_t$ adaptedness of $u_t^{l+1, M}$. In this work, we shall derive numerical methods to calculate the conditional estimated control process $\E[u_s^{l+1, M} | \cF^M_t]$ for $s \in [t, T]$, where $t$ gradually increases according to data reception. The procedure of calculating $\{\E[u_s^{l+1, M} | \cF^M_t]\}_{t \les s \les T}$ would give us an estimate for $u_t^{l+1, M}$, which is what we need in the gradient descent scheme \eqref{GD} to obtain $u^*_t$. However, in the conditional gradient descent iteration \eqref{E-GD_s}, the gradient $(J^{*})'_{u}$ on the right hand side is corresponding to $u_s^{l,M}$, which is $\cF^M_s$-adapted. Since the observational information $\{\cF^M_s\}_{t\les s \les T}$ is not available at time $t$, to make the iteration \eqref{E-GD_s} implementable, we use conditional expectation $\E[u_s^{l+1, M} | \cF^M_t]$ to replace $u_s^{l,M}$ since it provides the \textit{best approximation} for $u_s^{l+1, M}$ given $\cF^M_t$. In this way, we denote
$$u_s^{l, M}|_t := \E[u_s^{l+1, M} | \cF^M_t] $$
and carry out the following gradient descent iteration
\begin{equation}\label{GD_s}
u_s^{l+1, M}|_t = u_s^{l, M}|_t - \rho \E\big[ (J^{*})'_{u}(u_s^{l, M}|_t) | \cF^M_t\big], \ l = 0, 1, \cdots,   \ s \in [t, T],
\end{equation}
to derive $u_t^{l+1, M}$ as desired.
From equation \eqref{GD_s}, we can see that the conditional gradient process $\E\big[ (J^{*})'_{u}(u_s^{l,M}|_t) | \cF^M_t\big]$ is the key component and the effort to obtain such a term is composed by the following two tasks: (i) obtaining solutions $Y$ and $Z$ of the following FBSDEs system
\begin{equation}\label{FBSDE:system}
\begin{aligned}
 \ds dX_s&=b(s, X_s, u_s^{l, M}|_t)ds + \sigma(s, X_s, u_s^{l, M}|_t) dW_s, \qq \qq \qq \q s \in [t, T] \\
\ds dY_s &=\big(-b_x(s, X_s, u_s^{l, M}|_t)^\top Y_s-\si_x(s, X_s, u_s^{l, M}|_t)^\top Z_s\\
& \qq \qq \qq -f_x(s, X_s, u_s^{l, M}|_t)^\top \big)ds +Z_sdW_s + \z_s dB_s, \quad Y_T=h_x^\top(X_T)
\end{aligned}
\end{equation}
%
where the controlled process $X$ is corresponding to the estimated control $u_s^{l,M}|_t$; and (ii) obtaining effective evaluations for the conditional expectation $E[ \psi (s)| \cF^M_t]$ for a $\mathcal{F}_s^{X} \vee \mathcal{F}_t^{M}$ adapted stochastic process $\psi(s)$.
In most situations, the FBSDEs system \eqref{FBSDE:system} is not explicitly solvable. Therefore, we use numerical algorithms to approximate solutions $X$, $Y$ and $Z$ in this work.
The reason that we did not include the observation equation $dM_s=g(X_s)ds+dB_s$ in the FBSDEs system \eqref{FBSDE:system} is that $M$ is not contained in the expression \eqref{J'_t} for the gradient process. Moreover, when the control process is chosen to be $u_s^{l,M}|_t$, the adjoint BSDE in \eqref{FBSDE:system} is no longer driven by observational information beyond time $t$. In fact, $M$ is the data collected by observation equipments and the observational data $\{M_s\}_{s>t}$ is not available at time $t$, and therefore the observation dynamics should not be incorporated into the FBSDEs system.

On the other hand, when the controlled dynamics and the observation function $g$ are nonlinear, obtaining the conditional expectation $\E[\psi (s)|\cF^M_t]$ (in task (ii)) is not trivial either.
In this paper, we introduce an optimal filtering framework to achieve this goal.
The standard formulation for an optimal filtering problem is given by the following state-space model,
\begin{equation}\label{NLF}
\begin{aligned}
dS_t &= \beta(t, S_t) dt +\gamma(t, S_t) dW_t\\
d M_t & = g(S_t) dt + dB_t,
\end{aligned}
\end{equation}
where the drift term $\beta:\mathbb{R} \times \mathbb{R}^d \rightarrow \mathbb{R}^d$ describes some dynamical model, and $\gamma: \mathbb{R} \times \mathbb{R}^{d}\rightarrow \mathbb{R}^{d \times l}$ is the diffusion coefficient that drives the Brownian motion $W \in \mathbb{R}^l$ perturbing the model. The stochastic process $S_t$ is usually called the \textit{state process}, and $M_t$ provides partial noisy observations for $S_t$ in the same manner of the ``measurement process'' introduced in \eqref{observe[0,T]}. The goal of the optimal filtering problem is to obtain the best estimate for $\Psi(S_t)$ given the observational information $\cF^M_t$, where $\Psi$ is a test function. Mathematically, we aim to find the optimal filter $\tilde{\Psi}_t$ defined by the conditional expectation $\tilde{\Psi}_t : = \E[\Psi(S_t) | \cF^M_t]$. When all the functions $\beta$, $\gamma$ and $g$ in the optimal filtering problem are linear, the desired conditional expectation can be derived analytically by the well-known Kalman filter method. However, in most practical applications where the optimal filtering problems are nonlinear, numerical methods are used to approximate the conditional probability density function (pdf) for the state, i.e. $p(S_t | \cF^M_t)$, which can be used to calculate the optimal filter through integral
$\tilde{\Psi}_t = \int_{\mathbb{R}^d} \Psi(x) \cdot p(x \big| \cF^M_t) dx$. Well known nonlinear filtering methods include the Zakai filter, in which the desired filtering density is formulated as solution of a parabolic type stochastic PDE (Zakai equation), and the particle filter, which provides an empirical distribution to describe $p(S_t | \cF^M_t)$.

In our data driven feedback control problem, we let the controlled process $X_t$ be the target state process $S_t$ in the optimal filtering problem and generalize the test function $\Psi$ as a stochastic process defined by
\begin{equation}\label{test_filter}
\Psi(s, X_s, u_s^{l, M}|_t): = b_u(s, X_s, u_s^{l, M}|_t)^\top Y_s+\si_u(s, X_s, u_s^{l, M}|_t)^\top Z_s+f_u(s, X_s, u_s^{l, M}|_t)^\top
\end{equation}
for $s \in [t, T]$.
With the conditional pdf $p(X_t | \cF^M_t)$ that we obtain through optimal filtering methods and the fact that $\Psi(s, X_s, u_s^{l, M}|_t)$ is a stochastic process depending on the state of random variable $X_t$, the conditional gradient process $\E\big[ (J^{*})'_{u}(u_s^{l, M}|_t)| \cF^M_t\big]$ in \eqref{GD_s} can be obtained by the following integral
\begin{equation}\label{J':integral}
\E\big[ (J^{*})'_{u}(u_s^{l, M}|_t) | \cF^M_t\big] = \int_{\mathbb{R}^d} \E[\Psi(s, X_s, u_s^{l,M}|_t)\big|X_t=x] \cdot p(x \big| \cF^M_t) dx, \q s \in [t, T].
\end{equation}

In the following section, we derive a computational framework that carries out the gradient descent optimization \eqref{GD_s} to evaluate the data driven feedback control $u^*_t$, and introduce an efficient stochastic optimization algorithm to implement our computational framework. In many practical applications, it is difficult to control the diffusion. Therefore, in this work we assume that the diffusion coefficient $\sigma$ does not contain the control term.

\section{Numerical approach for data driven feedback control}\label{Numerical approach}

The general framework of our data driven feedback control system is designed on a temporal partition $\Pi_{N_T}$ that reflects the capability of collecting measurement data and implementing control actions, where $\Pi_{N_T}$ is defined by
$$\Pi_{N_T}: = \{t_n : 0 = t_0 < t_1 < t_2 < \cdots < t_n < \cdots < t_{N_T} = T\}, \q N_T \in \mathbb{N},$$
and we use the control sequence $\{u^{*}_{t_n}\}_{n=1}^{N_T}$ to represent the control process $u^*$.
In this way, the gradient descent \eqref{GD} is restricted to discrete time instants $\Pi_{N_T}$, and the conditional gradient descent scheme that we use to obtain $u^*_{t_n}$ becomes
\begin{equation}\label{GD_n_s}
u_{t_i}^{l+1, M}|_{t_n} = u_{t_i}^{l, M}|_{t_n} - \rho \E\big[ (J^{*})'_{u}(u_{t_i}^{l, M}|_{t_n}) | \cF^M_{t_n}\big],  \q n = 0, 1, \cdots, N_T, \q n \les i \les N_T,
\end{equation}
where the conditional gradient process on the right hand side of \eqref{GD_n_s} is given by
\begin{equation}\label{J'_tn}
\begin{aligned}
 \E[(J^*)'_u(u_{t_i}^{l,M}|_{t_n}) | \cF^M_{t_n}]   =  \E\Big[b_u(t_i, X_{t_i}, u_{t_i}^{l, M}|_{t_n})^\top Y_{t_i} +f_u(t_i, X_{t_i}, u_{t_i}^{l, M}|_{t_n})^\top\bigm|\cF^M_{t_n}\Big],
\end{aligned}
\end{equation}
and $X_{t_i}$ is $\cF^W_{t_i}\vee \cF^M_{t_n}$ measurable since its corresponding control is $u_{t_i}^{l, M}|_{t_n}$.
In what follows, we introduce numerical schemes to calculate $X_{t_i}$, $Y_{t_i}$ and $Z_{t_i}$ corresponding to the estimated conditional optimal control $u_{t_i}^{l, M}|_{t_n}$ in Section \ref{Numerical_FBSDEs}, and introduce a particle filter method, which is one of the most widely accepted methods for solving nonlinear filtering problems, to approximate the conditional distribution $p(X_{t_n}|\cF^M_{t_n})$ in Section \ref{PF_distribution}. Then, in Section \ref{PF-SGD} we combine schemes in Sections \ref{Numerical_FBSDEs} - \ref{PF_distribution} to formulate a computational framework for the data driven feedback control problem with an efficient stochastic optimization implementation, and we summarize our numerical approach in Section \ref{Summary}.

\subsection{Numerical schemes for FBSDEs}\label{Numerical_FBSDEs}

We consider the FBSDEs system \eqref{FBSDE:system} on time interval $[t_i, t_{i+1}]$, i.e.
\begin{equation}\label{FBSDEs_ti}
\begin{aligned}
 X_{t_{i+1}}=& X_{t_i} + \int_{t_i}^{t_{i+1}}b(s,X_s,u_s^{l,M}|_{t_n})ds+\int_{t_i}^{t_{i+1}}\si(s,X_s)dW_s,  \qq \q\text{(SDE)}\\
Y_{t_i}=& Y_{t_{i+1}} + \int_{t_i}^{t_{i+1}}\Big(b_x(s,X_s,u_s^{l,M}|_{t_n})^\top Y_s+\si_x(s,X_s)^\top Z_s+f_x(s,X_s,u_s^{l,M}|_{t_n})^\top\Big)ds \\
 & \qq \qq -\int_{t_i}^{t_{i+1}} Z_sdW_s + \int_{t_i}^{t_{i+1}} \z_sdB_s.  \qq \qq \qq \qq \ \text{(BSDE)}
\end{aligned}
\end{equation}

To solve the first equation in \eqref{FBSDEs_ti}, which is a standard forward SDE, we use the left-point formula to approximate both the deterministic and stochastic integrals, and obtain the following Euler-Maruyama discretized approximation equation for $X$ as
\begin{equation}\label{Dis:X}
X_{t_{i+1}}=X_{t_i} + b(t_i,X_{t_i},u_{t_i}^{l,M}|_{t_n}) \Delta t_i+\si(t_i,X_{t_i})\Delta W_{t_i} + R_{X, n}^i,
\end{equation}
where $R_{X,n}^i$ is the approximation error for integrals, and we denote $\Delta t_i := t_{i+1} - t_i$ and $\Delta W_{t_i} := W_{t_{i+1}} - W_{t_i}$.

To derive a discretization scheme for solution $Y$, we take conditional expectation $\E_i[\cdot]:=\E[\cdot|  \cF^{X,B}_{t_i}]$ on both sides of the BSDE in \eqref{FBSDEs_ti} to get
$$Y_{t_i}= \E_i[Y_{t_{i+1}}] + \int_{t_i}^{t_{i+1}}\E_i\big[ b_x(s,X_s,u_s^{l,M}|_{t_n})^\top Y_s+\si_x(s,X_s)^\top Z_s+f_x(s,X_s,u_s^{l,M}|_{t_n})^\top \big]ds, $$
where the stochastic integrals are eliminated and we have used the fact that $Y_{t_i} = \E_i[Y_{t_i}]$ due to the adaptedness of $Y$. Then, we use the right-point formula to approximate the deterministic integral in the above equation and obtain
\begin{equation}\label{Exp:FBSDEs_ti}
\begin{aligned}
Y_{t_i}= \E_i[Y_{t_{i+1}}] + & \E_i\Big[ b_x(t_{i+1},X_{t_{i+1}},u_{t_{i+1}}^{l,M}|_{t_n})^\top Y_{t_{i+1}}+\si_x(t_{i+1},X_{t_{i+1}})^\top Z_{t_{i+1}} \\
 & \qq +f_x(t_{i+1},X_{t_{i+1}},u_{t_{i+1}}^{l,M}|_{t_n})^\top \Big] \Delta t_i + R_{Y, n}^i,
\end{aligned}
\end{equation}
where $R_{Y,n}^i$ is the approximation error for the integral.

For solution $Z$, we multiply $\Delta W_{t_i}$ and take conditional expectation $\E_i[\cdot]$ on both sides of the BSDE in  \eqref{FBSDEs_ti} to get
\begin{equation*}\label{Exp:WZ}
\begin{aligned}
\E_i\big[\int_{t_i}^{t_{i+1}}Z_s  \Delta W_{t_i}  dW_s \big]  =& \E_{i}[Y_{t_{i+1}} \Delta W_{t_i}] + \int_{t_i}^{t_{i+1}}\E_i\Big[ \Big( b_x(s,X_s,u_s^{l,M}|_{t_n})^\top Y_s \\
& \ +\si_x(s,X_s)^\top Z_s+f_x(s,X_s,u_s)^\top \Big)\Delta W_{t_i}\Big] ds,
\end{aligned}
\end{equation*}
where we have used the fact $ \E_{i}[Y_{t_{i}} \Delta W_{t_i}] = 0$ and the independence of Brownian motions $W$ and $B$. Then, we use the left-point formula to approximate both deterministic and stochastic integrals in the above equation and obtain
\begin{equation}\label{Exp:FBSDEs-W_ti}
\begin{aligned}
Z_{t_i}  \Delta t_{i} =& \ \E_{i}[Y_{t_{i+1}} \Delta W_{t_i}] + R_{Z,n}^i,
\end{aligned}
\end{equation}
where the deterministic integral is eliminated due to the $\dbF^W$ adaptedness and the fact that we picked the left-point formula approximation, and $R_{Z,n}^i$ is the term containing approximation errors for integrals.

By dropping approximation error terms $R_{X,n}^i$, $R_{Y,n}^i$ and $R_{Z,n}^i$ in \eqref{Dis:X}, \eqref{Exp:FBSDEs_ti} and \eqref{Exp:FBSDEs-W_ti}, respectively, we could provide a numerical algorithm for solving $X$, $Y$ and $Z$ corresponding to the gradient process $\E[(J^*)'_u(u_{t_i}^{l,M}|_{t_n}) | \cF^M_{t_n}]$ as introduced in \eqref{J'_tn}:\\

\noindent For a fixed $n\les N_T-1$ representing the current time $t_n$, and for $i = N_{T}-1, \cdots, n+1, n$ with initial conditions $Y_{N_T}$ and $Z_{N_T}$, we solve the FBSDEs system through the following schemes
\begin{equation}\label{Semi-scheme}
\begin{aligned}
X_{i+1} =& X_i + b(t_i, X_i, u_{t_i}^{l,M}|_{t_n}) \Delta t_i + \sigma(t_i, X_i) \Delta W_{t_i}\\
Y_{i} = & \E_i[Y_{i+1}] +  \E_i\Big[ b_x(t_{i+1},X_{i+1},u_{t_{i+1}}^{l,M}|_{t_n})^\top Y_{i+1}+\si_x(t_{i+1}, X_{i+1})^\top Z_{i+1} \\
 & \qq +f_x(t_{i+1},X_{i+1}, u_{t_{i+1}}^{l,M}|_{t_n})^\top \Big] \Delta t_i \\
 Z_{i}  =& \ \E_{i}[Y_{i+1} \Delta W_{t_i}] \cdot (\Delta t_{i})^{-1}
 \end{aligned}
\end{equation}
where $X_{i+1}$, $Y_i$ and $Z_i$ are numerical approximations for $X_{t_{i+1}}$, $Y_{t_i}$ and $Z_{t_i}$, respectively. The numerical analysis for schemes \eqref{Semi-scheme} and their extended versions are discussed in \cite{Bao_first, Zhao_multi}.

\begin{rem}
We notice that the random variable $\z$ is not involved in numerical schemes \eqref{Semi-scheme} to calculate approximate solutions $X_{i+1}$, $Y_i$ and $Z_i$, and it's not contained in the expression for $(J^*)'_u$. Therefore, we do not need an extra scheme to approximate $\z$ in our approach.
\end{rem}
\vspace{0.25em}


In order to carry out numerical schemes \eqref{Semi-scheme}, one needs to evaluate conditional expectations. There are many well-known methods to approximate expectations, such like (quasi) Monte Carlo methods and numerical integration methods, etc. In this work, we adopt the standard Monte Carlo method due to its high effectiveness and efficiency in approximating expectations -- especially in high dimensional spaces. Specifically, we use a set of $K \in \mathbb{N}$ random samples  to describe $\Delta W_{t_i}$ and the expectations are approximated by sample means. As a result, schemes \eqref{Semi-scheme} now become
\begin{equation}\label{MC-scheme}
\begin{aligned}
X^k_{i+1} =& X_i + b(t_i, X_i, u_{t_i}^{l,M}|_{t_n}) \Delta t_i + \sigma(t_i, X_i) \sqrt{\Delta t_i} \omega^k_i, \qq  k = 1,2 \cdots, K, \\
Y_{i} = & \sum_{k=1}^{K}\frac{Y_{i+1}^k}{K} +  \frac{\Delta t_i}{K} \sum_{k=1}^{K} \Big[ b_x(t_{i+1},X^k_{i+1},u_{t_{i+1}}^{l,M}|_{t_n})^\top Y^k_{i+1} \\
 & \qq \q +\si_x(t_{i+1}, X^k_{i+1})^\top Z^k_{i+1}   +f_x(t_{i+1},X^k_{i+1}, u_{t_{i+1}}^{l,M}|_{t_n})^\top \Big], \\
 Z_{t_i}   =& \ \frac{1}{\Delta t_{i}}\sum_{k=1}^{K}\frac{Y_{i+1}^k \sqrt{\Delta t_i} \omega^k_i}{K},
 \end{aligned}
\end{equation}
where $\{\omega^k_i\}_{k=1}^K$ is a set of random samples following the standard Gaussian distribution that we use to describe the randomness of $\Delta W_{t_i}$; $Y^k_{i+1}$ and $Z^k_{i+1}$ denote approximate solutions $Y_{i+1}$ and $Z_{i+1}$ corresponding to the $k$-th random sample in the Monte Carlo approximation.

%
%

\subsection{Particle filter method for conditional distribution}\label{PF_distribution}
The purpose of our optimal filtering procedure is to obtain the conditional distribution for the controlled process, which will be used to generate the conditional gradient process $\E\big[ (J^{*})'_{u}(u_{t_i}^{l, M}|_{t_n}) | \cF^M_{t_n}\big]$ with respect to the control process. To proceed, we consider the controlled process on time interval $[t_{n-1}, t_n]$, i.e.
\begin{equation}\label{NLF-state}
X_{t_{n}} = X_{t_{n-1}} + \int_{t_{n-1}}^{t_n}b(s, X_s, u_s) ds + \int_{t_{n-1}}^{t_{n}}\sigma(s, X_s) d W_s,
\end{equation}
and the observational data at time instant $t_{n}$ is given by $M_{t_{n}}  = g(X_{t_{n}}) + \eta_{n}$,
where $\eta_n$ is the observational noise.
Well known pioneer work by Zakai showed that the desired filtering density is equivalent to the solution of a stochastic PDE (the Zakai equation). Although the Zakai's approach theoretically solves the nonlinear filtering problem, due to the ``curse of dimensionality'' of solving a stochastic PDE in high dimensional space, the Zakai filter is prohibitive in practically applications.
In this work, we adopt the Bayesian filter framework and introduce a widely accepted Bayesian filter -- the particle filter method -- to solve for the filtering density  $p(X_{t_{n}} | \cF^M_{t_n})$. The Bayesian filter framework is composed by two stages: \textit{prediction stage} and \textit{update stage}. In the prediction stage, assuming that we have the distribution $p(X_{t_{n-1}} | \cF^M_{t_{n-1}})$ at time instant $t_{n-1}$, the prior pdf that predicts the state of controlled process at time $t_{n}$ is given by the following Chapman-Kolmogorov formula
\begin{equation*}
p(X_{t_{n}} | \cF^M_{t_{n-1}}) = \int p(X_{t_{n-1}} | \cF^M_{t_{n-1}} ) p(X_{t_{n}} | X_{t_{n-1}}) dX_{t_{n-1}},
\end{equation*}
where $p(X_{t_{n}} | X_{t_{n-1}})$ is the transition probability derived from the state dynamics \eqref{NLF-state}. With the new observational data $M_{t_{n}}$, the update stage uses the Bayesian inference to update the prior pdf and get the posterior pdf $p(X_{t_{n}} | \cF^M_{t_{n}})$ as
\begin{equation}\label{Bayesian}
p(X_{t_{n}} | \cF^M_{t_{n}}) = \frac{p(X_{t_{n}} | \cF^M_{t_{n-1}}) p(M_{t_{n}} | X_{t_{n}})}{p(M_{t_{n}} | \cF^M_{t_{n-1}})},
\end{equation}
where $p(M_{t_{n}} | X_{t_{n}})$ is the likelihood function that describes the discrepancy between the predicted state and the observations \cite{Do2}. In what follows, we introduce the benchmark particle filter method (bootstrap filter algorithm \cite{particle-filter}) to the implement the aforementioned Bayesian filter framework due to the high effectiveness and efficiency of the particle filter in solving nonlinear filtering problems.  To proceed, assume that at time instant $t_{n-1}$ we have $S$ particles, denoted by $\{x_{n-1}^{(s)}\}_{s=1}^S$, that follow an empirical distribution $\pi(X_{t_{n-1}}|\cF^M_{t_{n-1}}):= \frac{1}{S} \sum_{s=1}^S \delta_{x_{n-1}^{(s)}}(X_{t_{n-1}})$ as an approximation for $p(X_{t_{n-1}} | \cF^M_{t_{n-1}})$, where $\delta_x$ is the Dirac delta function at $x$. The prior pdf that we want to find in the prediction stage is approximated as
\begin{equation}\label{PF-predict}
\tilde{\pi}(X_{t_{n}} | \cF^M_{t_{n-1}}) := \frac{1}{S} \sum_{s=1}^S \delta_{\tilde{x}_{n}^{(s)}} (X_{t_{n}}),
\end{equation}
where $\tilde{x}_{n}^{(s)}$ is sampled from $\pi(X_{t_{n-1}} | \cF^M_{t_{n-1}}) p(X_{t_{n}} | X_{t_{n-1}})$. As a result, the sample cloud $\{\tilde{x}_{n}^{(s)}\}_{s=1}^S$ provides an approximate distribution for the prior $p(X_{t_{n}} | \cF^M_{t_{n-1}})$. Then, in the update stage, we update the approximated prior pdf to get the posterior pdf by replacing $p(X_{t_{n}} | \cF^M_{t_{n-1}})$ with $\tilde{\pi}(X_{t_{n}} | \cF^M_{t_{n-1}})$ in the Bayesian inference \eqref{Bayesian} as following
\begin{equation}
\begin{aligned}
\tilde{\pi}(X_{t_{n}} | \cF^M_{t_{n}}) := \frac{\sum_{s=1}^S \delta_{\tilde{x}_{n}^{(s)}}(X_{t_{n}}) p(M_{t_{n}} | \tilde{x}_{n}^{(s)})}{\sum_{s=1}^S p(M_{t_{n}} | \tilde{x}_{n}^{(s)})} =  \sum_{s=1}^S w_{n}^{(s)} \delta_{\tilde{x}_{n}^{(s)}}(X_{t_{n}}).
\end{aligned}
\end{equation}
In this way, we obtain a weighted empirical distribution $\tilde{\pi}(X_{t_{n}} | \cF^M_{t_{n}})$ that approximates the posterior pdf $p(X_{t_{n}} | \cF^M_{t_{n}})$ with the importance density weight $w_{n}^{(s)} \propto p(M_{t_{n}} | \tilde{x}_{n}^{(s)})$. In practice, importance weights $\{w_{n}^{(s)}\}_{s=1}^S$ tend to concentrate on a few samples after several time steps, which dramatically reduces the effective particle size in the algorithm.
To avoid such degeneracy problem, we resample the particles $\{\tilde{x}_{n}^{(s)}\}_{s=1}^S$ by replacing particles with low density weights with copies of particles with high weights. In the bootstrap particle filter, we use importance sampling method to generate equally weighted samples  $\{x_{n}^{(s)}\}_{s=1}^S$ from $\tilde{\pi}(X_{t_{n}} | \cF^M_{t_n})$ to formulate our empirical distribution $\pi(X_{t_{n}} | \cF^M_{t_n})$ that describes the conditional pdf for the controlled process $p(X_{t_{n}} | \cF^M_{t_n})$ \cite{particle-filter}, i.e.
\begin{equation}\label{PF-resample}
\pi(X_{t_{n}} | \cF^M_{t_n})= \frac{1}{S} \sum_{s=1}^S \delta_{x_{n}^{(s)}}(X_{t_{n}}).
\end{equation}


\subsection{Stochastic optimization for control process}\label{PF-SGD}
The procedure that aims to determine the data driven optimal control is the gradient descent iteration. In this subsection, we combine the numerical schemes for the adjoint FBSDEs system \eqref{FBSDE:system} (introduced in Section \ref{Numerical_FBSDEs}) and the particle filter algorithm (introduced in Section \ref{PF_distribution}) to formulate an efficient stochastic optimization algorithm to solve for the optimal control process $u^*$.

The framework that we adopt to design our algorithm is the conditional gradient descent iteration \eqref{GD_n_s} and the conditional gradient process with respect to the control process, i.e. $\E[(J^*)'_u(u_{t_i}^{l,M}|_{t_n}) | \cF^M_{t_n}]$, is calculated by using the integral expression \eqref{J':integral}. Specifically, on a time instant $ t_n \in \Pi_{N_T}$, we have
\begin{equation*}\label{J'_tn:integral}
\E\big[ (J^{*})'_{u}(u_{t_i}^{l, M}|_{t_n}) | \cF^M_{t_n}\big] = \int_{\mathbb{R}^d} \E[\Psi(t_i, X_{t_i}, u_{t_i}^{l,M}|_{t_n})\big|X_{t_n}=x] \cdot p(x \big| \cF^M_{t_n}) dx,
\end{equation*}
where $t_i \geq t_n$ is a time instant after $t_n$ and the test function $\Psi$ in \eqref{test_filter} is chosen as
$$\Psi(t_i, X_{t_i}, u_{t_i}^{l,M}|_{t_n}) = b_u(t_i, X_{t_i}, u_{t_i}^{l, M}|_{t_n})^\top Y_{t_i}+f_u(t_i, X_{t_i}, u_{t_i}^{l, M}|_{t_n})^\top.$$
In our numerical approach, we use approximate solutions $(Y_i, Z_i)$ of FBSDEs from schemes \eqref{MC-scheme} to replace $(Y_{t_i}, Z_{t_i})$ and the conditional distribution $p(X_{t_n} \big| \cF^M_{t_n})$ is approximated by the empirical distribution  $\pi(X_{t_{n}} | \cF^M_{t_n})$ obtained from the particle filter algorithm \eqref{PF-predict} - \eqref{PF-resample}. Specifically, for increasing time instants $t_n \in \Pi_{N_T}$, $n = 0, 1, 2, \cdots, N_T$, we solve for the optimal control $u^*_{t_n}$ through the following gradient descent optimization iteration
\begin{equation}\label{GD_tn}
u_{t_i}^{l+1, M}|_{t_n} = u_{t_i}^{l, M}|_{t_n} - \rho \tilde{\E}_{t_n}\big[ (J^{*})'_{u}(u_{t_i}^{l, M}|_{t_n})\big],  \q l = 0, 1, 2, \cdots, L, \q t_i \geq t_n,
\end{equation}
where $\tilde{\E}_{t_n}\big[ (J^{*})'_{u}(u_{t_i}^{l, M}|_{t_n})\big]$ is defined by
\begin{equation}\label{Approx-J_tn}
\begin{aligned}
 \tilde{\E}_{t_n}\big[ (J^{*})'_{u}(u_{t_i}^{l, M}|_{t_n})\big] := & \frac{1}{S}\sum_{s = 1}^{S}\E\Big[b_u(t_i, X_{t_i}, u_{t_i}^{l, M}|_{t_n})^\top Y_{i} \\
 & \qq \qq +f_u(t_i, X_{t_i}, u_{t_i}^{l, M}|_{t_n})^\top \bigm| X_{t_n} = x_n^{(s)}\Big],
\end{aligned}
\end{equation}
and samples $\{x_n^{(s)}\}_{s=1}^S$ follow distribution $\pi(X_{t_{n}} | \cF^M_{t_n})$ which describes the conditional pdf of the controlled process $X_{t_n}$ given the observational information $\cF^M_{t_n}$. Here, we want to re-emphasize that the purpose of carrying out an optimization procedure for the conditional control process $\{u_{t_i}^{\cdot, M}|_{t_n}\}_{i=n}^{N_T}$ is to find an estimate for the optimal control at time $t_n$. Since the approximate conditional gradient process on the right hand side of \eqref{Approx-J_tn} is still under expectation, Monte Carlo simulation for $\E[\cdot |  X_{t_n} = x_n^{(s)}]$ is needed in general, and $\tilde{\E}_{t_n}\big[ (J^{*})'_{u}(u_{t_i}^{l, M}|_{t_n})\big]$ is approximated as following
\begin{equation}\label{MC-J_tn}
\begin{aligned}
 \tilde{\E}_{t_n}\big[ (J^{*})'_{u}(u_{t_i}^{l, M}|_{t_n})\big] \approx & \frac{1}{S} \frac{1}{\Lambda} \sum_{s = 1}^{S} \sum_{\lambda=1}^{\Lambda} \Big[b_u(t_i, X^{(\lambda,s)}_{t_i}, u_{t_i}^{l, M}|_{t_n})^\top Y^{(\lambda,s)}_{i} \\
 & \qq \qq +f_u(t_i, X^{(\lambda, s)}_{t_i}, u_{t_i}^{l, M}|_{t_n})^\top \bigm| X_{t_n} = x_n^{(s)}\Big],
\end{aligned}
\end{equation}
where $\Lambda$ is the number of samples to approximate the expectation $\E[\cdot |  X_{t_n} = x_n^{(s)}]$, $X^{(\lambda,s)}_{t_i}$ is the $\lambda$-th realization of controlled process with initial state $X_{t_n} = x_n^{(s)}$, and $Y_i^{(\lambda,s)}$ is the approximate solution $Y_i$ corresponding to $X^{(\lambda,s)}_{t_i}$. We can see from the above Monte Carlo approximation that in order to approximate the expectation $\tilde{\E}_{t_n}\big[ (J^{*})'_{u}(u_{t_i}^{l, M}|_{t_n})\big]$ in one gradient descent iteration step, we need to generate altogether $S \times \Lambda$ samples of controlled process and evaluate the values of solution $Y_i$ corresponding to $\{X^{(\lambda,s)}_{t_i}\}_{\lambda, s}$ as well. This is even more computationally expensive when the controlled system is a high dimensional process.

Inspired by the application of stochastic gradient descent (SGD) in improving the efficiency of classic gradient descent optimization, in this work we introduce a stochastic optimization algorithm that combines the particle filter method with the stochastic gradient descent algorithm to carry out the gradient descent iteration procedure \eqref{GD_tn} efficiently. When implementing the stochastic gradient descent in the optimization scheme \eqref{GD_tn}, instead of using the fully calculated Monte Carlo simulation to approximate the conditional expectation on the right hand side of \eqref{Approx-J_tn} (as introduced in \eqref{MC-J_tn}), we use only one realization of $X_{t_i}$ to represent the expectation, i.e.
\begin{equation}\label{SGD-J_tn}
\begin{aligned}
 & \E\Big[b_u(t_i, X_{t_i}, u_{t_i}^{l, M}|_{t_n})^\top Y_{i}  +f_u(t_i, X_{t_i}, u_{t_i}^{l, M}|_{t_n})^\top \bigm| X_{t_n} = x_n^{(s)}\Big] \\
 \approx & \ b_u(t_i, X^{(\hat{l},s)}_{t_i}, u_{t_i}^{l, M}|_{t_n})^\top Y^{(\hat{l},s)}_{i}+f_u(t_i, X^{(\hat{l},s)}_{t_i}, u_{t_i}^{l, M}|_{t_n})^\top,
\end{aligned}
\end{equation}
where $X^{(\hat{l},s)}_{t_i}$ is a randomly generated realization of the controlled process with initial state $X^{(\hat{l},s)}_{t_n} = x_n^{(s)}$, and we use the index $\hat{l}$ to indicate that the random generation of the controlled process varies among the gradient descent iteration steps in \eqref{GD_tn}. In addition, one should notice that the particles $\{ x_n^{(s)}\}_{s=1}^S$ are used to describe the conditional distribution of the controlled process and the Monte Carlo average in the the approximation scheme \eqref{Approx-J_tn} also aims to approximate a conditional expectation. In light of the concept of single-realization representation for expectation as we described in \eqref{SGD-J_tn}, we extend this concept to the Monte Carlo averaging of the particles in the approximation scheme \eqref{Approx-J_tn}. Therefore, we use the following expresssion to represent the conditional expectation $\tilde{\E}_{t_n}\big[ (J^{*})'_{u}(u_{t_i}^{l, M}|_{t_n})\big]$ in \eqref{GD_tn}
\begin{equation}\label{Exp-SGD}
\begin{aligned}
\tilde{\E}_{t_n}\big[ (J^{*})'_{u}(u_{t_i}^{l, M}|_{t_n})\big] \approx b_u(t_i, X^{(\hat{l},\hat{s})}_{t_i}, u_{t_i}^{l, M}|_{t_n})^\top Y^{(\hat{l},\hat{s})}_{i}+f_u(t_i, X^{(\hat{l},\hat{s})}_{t_i}, u_{t_i}^{l, M}|_{t_n})^\top,
\end{aligned}
\end{equation}
where $X^{(\hat{l},\hat{s})}_{t_i}$ indicates a \textit{randomly generated} realization of the controlled process with a \textit{randomly selected} initial state $X^{(\hat{l},\hat{s})}_{t_n} = x_n^{(\hat{s})}$ from the particle cloud $\{x_n^{(s)}\}_{s=1}^S$.
In this way, the gradient descent optimization scheme \eqref{GD_tn} becomes the following SGD optimization iteration scheme

\begin{equation}\label{SGD_tn}
\begin{aligned}
u_{t_i}^{l+1, M}|_{t_n} = & u_{t_i}^{l, M}|_{t_n} -  \rho \Big(b_u(t_i, X^{(\hat{l},\hat{s})}_{t_i}, u_{t_i}^{l, M}|_{t_n})^\top Y^{(\hat{l},\hat{s})}_{i} \\
& \qq \qq + f_u(t_i, X^{(\hat{l},\hat{s})}_{t_i}, u_{t_i}^{l, M}|_{t_n})^\top \Big),   \q l = 0, 1, 2, \cdots, L, \q t_i \geq t_n,
\end{aligned}
\end{equation}
where $Y^{(\hat{l},\hat{s})}_{i}$ is the approximate solution $Y_i$ corresponding to the random sample $X^{(\hat{l},\hat{s})}_{t_i}$ indexed by the sample pair $(\hat{l},\hat{s})$, and $X^{(\hat{l},\hat{s})}_{t_i}$ is generated as following
\begin{equation}\label{Sample:X}
X^{(\hat{l},\hat{s})}_{t_{i+1}} = X^{(\hat{l},\hat{s})}_{t_i} + b(t_i, X^{(\hat{l},\hat{s})}_{t_i}, u_{t_i}^{l, M}|_{t_n}) \Delta t_i + \sigma(s, X_s) \sqrt{\Delta t_i} \ \omega^{(\hat{l},\hat{s})}_i, \ i = n, \cdots, N_T-1
\end{equation}
where $X^{(\hat{l},\hat{s})}_{t_n} = x_n^{(\hat{s})} \in \{x_n^{(s)}\}_{s=1}^S$, $\omega^{(\hat{l},\hat{s})}_i \sim N(0, 1)$, and $\{\omega^{(\hat{l},\hat{s})}_i\}_{i=1}^{N_T-1}$ form a sequence of random variables corresponding to the sample indices $\hat{l}$ and $x_n^{(\hat{s})}$. Then, the resulted conditional optimal control process $\{u_{t_i}^{L, M}|_{t_n}\}_{i=n}^{N_T}$ leads to an estimate for our desired data driven optimal control at time instant $t_n$, i.e. $\hat{u}^*_{t_n}$, which is given by
$$\hat{u}^*_{t_n}:=  u_{t_n}^{L, M}|_{t_n}.$$
As a result, particles in the particle filter method (PF) and random samples in the stochastic gradient descent method (SGD) are combined into one stochastic optimization iteration \eqref{SGD_tn},
and we use only one-realization to represent two layers of expectations instead of altogether $S \times \Lambda$ calculations of $X^{(\lambda, s)}_{t_i}$ and $Y^{(\lambda, s)}_{i}$. For convenience of presentation, we name this methodology the ``PF-SGD''.

It is important to point out that in order to characterize the dependence of numerical solution $Y_{i}$ on $X_{t_i}^{(\hat{l},\hat{s})}$ as we indicated in \eqref{SGD_tn}, one needs to provide approximations for $Y$ and $Z$ corresponding to the random variable $X$. In practice, people typically derive approximate solutions $Y_i$ and $Z_i$ on a set of spatial points that represent $X_{t_i}$, and then construct interpolatory approximations based on discrete solution values on those selected spatial points. This may also be considered as spatial dimension approximation for numerical solutions of FBSDEs. Recently, a machine learning based method has been developed, which could solve FBSDEs in high dimensional space on a set of pre-selected spatial points \cite{ML_BSDE}. However, even with accurately calculated solutions on a large number of spatial points, characterizing $Y$ and $Z$ with respect to the random variable $X$ is still challenging due to the sparsity of spatial points in high dimensional space. To address the challenge in spatial dimension approximation, in this work we take the advantage of SDE nature of FBSDEs and incorporate numerical schemes for FBSDEs into our PF-SGD framework. Specifically, we further extend the concept of singel-realization representation of expectation to the numerical schemes \eqref{MC-scheme} for FBSDEs and use the solution path $\{X^{(\hat{l},\hat{s})}_{t_i}\}_{i=n}^{N_T}$ to drive our numerical solution for $Y$ and $Z$ as following
\begin{equation}\label{SGD-scheme}
\begin{aligned}
Y^{(\hat{l},\hat{s})}_{i} = & Y_{i+1}^{(\hat{l},\hat{s})} +  \Big[ b_x(t_{i+1},X^{(\hat{l},\hat{s})}_{t_{i+1}},u_{t_{i+1}}^{l,M}|_{t_n})^\top Y^{(\hat{l},\hat{s})}_{i+1} \\
 & \qq +\si_x(t_{i+1}, X^{(\hat{l},\hat{s})}_{t_{i+1}})^\top Z^{(\hat{l},\hat{s})}_{i+1}   +f_x(t_{i+1},X^{(\hat{l},\hat{s})}_{i+1}, u_{t_{i+1}}^{l,M}|_{t_n})^\top \Big] \Delta t_i, \\
 Z^{(\hat{l},\hat{s})}_{i} =& \ Y^{(\hat{l},\hat{s})}_{i+1} \omega^{(\hat{l},\hat{s})}_i  ( \Delta t_i)^{-1/2} ,
 \end{aligned}
\end{equation}
where the Monte Carlo type approximations for expectations are replaced by a single sample path corresponding to
$X^{(\hat{l},\hat{s})}_{t_i}$.
Although the schemes \eqref{SGD-scheme} does not provide accurate numerical approximations for solutions of FBSDEs -- compared with the conventional fully calculated schemes, such like the schemes \eqref{MC-scheme} and schemes in \cite{Touzi_BSDE, Four_step, Zhao_multi}, we want to point out that the conditional expectations in numerical schemes \eqref{Semi-scheme} for $Y_i$ and $Z_i$ appear in the expression for the gradient process $\tilde{\E}_{t_n}\big[ (J^{*})'_{u}(u_{t_i}^{l, M}|_{t_n})\big]$ in \eqref{GD_tn}. Therefore, the justification for the application of SGD in the gradient descent method also applies to \eqref{SGD-scheme}. In other words, the purpose of obtaining numerical solutions $Y_i$ and $Z_i$ for the adjoint FBSDEs system is to carry out the gradient descent optimization -- not to derive accurate approximations for FBSDEs.
In this way, the PF-SGD schemes \eqref{SGD_tn} - \eqref{SGD-scheme} combined together transfer our computational efforts in solving FBSDEs to optimization iterations of finding the optimal control.

\subsection{Summary of the numerical approach}\label{Summary}
In \textit{Algorithm 1} (Table \ref{Algorithm}), we summarize our PF-SGD algorithm for the data driven feedback control problem (C$^*$).
\vspace{0.5em}

\begin{table}\caption{}\label{Algorithm}
\vspace{1em}
\centering
\begin{tabular} {p{0.9\textwidth}}
\hline\noalign{\smallskip}
{\bf Algorithm 1}: {\em PF-SGD algorithm for data driven feedback control problem }\\
\noalign
{\smallskip}\hline
\noalign{\smallskip}
\vspace{-0.1cm}
\begin{spacing}{1.1}
\begin{algorithmic}\label{algorithm}
\item[Initialize] the particle cloud $\{x_0^{(s)}\}_{s=1}^S \sim \xi $ and the number of iteration $L \in \mathbb{N}$
\item[\textbf{while}] $n =0, 1, 2, \cdots, N_T$, \textbf{do} \\
\begin{description}
\item \hspace{-0.75em}  Initialize an estimated control process $\{u_{t_i}^{0, M}|_{t_n}\}_{i=n}^{N_T}$ and a step-size $\rho$;
\item[\textbf{for}] SGD iteration steps $l = 0, 1, 2, \cdots, L$,   \\
\begin{description}
\item \hspace{-2.5em}  Simulate one realization of controlled process $\{X^{(\hat{l},\hat{s})}_{t_{i+1}}\}_{i=n}^{N_T-1}$ through scheme \eqref{Sample:X} with $X^{(\hat{l},\hat{s})}_{t_n} = x_n^{(\hat{s})} \in \mathcal \{x_n^{(s)}\}_{s=1}^S$ ;
\item \hspace{-2.5em}  Calculate solution  $\{Y^{(\hat{l},\hat{s})}_{i}\}_{i=N_T}^{n}$ of the FBSDEs system \eqref{FBSDE:system} corresponding to $\{X^{(\hat{l},\hat{s})}_{t_{i+1}}\}_{i=n}^{N_T-1}$ through schemes \eqref{SGD-scheme};
\item \hspace{-2.5em}  Update the control process to obtain $\{u_{t_i}^{l+1, M}|_{t_n}\}_{i=n}^{N_T}$ through scheme \eqref{SGD_tn};
\end{description}
\item[\textbf{end for}]
\item \hspace{-0.75em} The estimated optimal control is given by $\hat{u}^*_{t_n} =  u_{t_n}^{L, M}|_{t_n}$;
\item \hspace{-0.75em}  Propagate particles through the particle filter algorithm \eqref{PF-predict} - \eqref{PF-resample} to obtain $\{x_{n+1}^{(s)}\}_{s=1}^S$ by using the estimated optimal control $\hat{u}^*_{t_n}$.
\end{description}
\item[\textbf{end while}]
\end{algorithmic}
\vspace{-1.2em}
\end{spacing}\\
\hline
\end{tabular}
\end{table}
The methodology of our approach is based on the formulation of a conditional optimization \eqref{GD_s}, which provides the ``best approximation'' for the data driven feedback control $u^*$ given the current observational information. The general computational framework combines (a) numerical schemes for FBSDEs, which provide numerical approximation for the gradient process with respect to the optimal control, and (b) the optimal filtering algorithm for the conditional distribution of the controlled process, which is used to computing the conditional expectation of the gradient process. Then, we adopt the concept of SGD and use a single-realization of simulated controlled process to represent expectations in our computational framework. The contribution of the PF-SGD algorithm is that every realization of particle based simulated controlled process $X^{(\hat{l},\hat{s})}_{t_{i}}$, as well as its corresponding FBSDEs solution pair $( Y^{(\hat{l},\hat{s})}_{i}, Z^{(\hat{l},\hat{s})}_{i} )$, is effectively used to search for the optimal control, which significantly saves unnecessary computational cost on solving the optimal control problem in the entire controlled state space.

\section{Numerical Experiments}\label{Numerics}

In this section, we present three numerical examples to demonstrate the effectiveness and efficiency of our data driven feedback control algorithm.  The processor that we use to run the numerical experiments is an Intel Core i7 with 2.5 GHz Dual-Core. In the first example, we demonstrate the effectiveness of our computational framework by solving a classic linear-quadratic (LQ) optimal control problem with \textit{nonlinear} observations. If the real state of controlled process is given, its corresponding optimal control process can be derived analytically. Therefore, this example could be a benchmark test that compares the data driven optimal control calculated from our algorithm with the analytical optimal control process, and we want to use this example to show that our algorithm could accurately provide the optimal control process based on indirect observations of the controlled process.  In the second example, we demonstrate the efficiency of the PF-SGD implementation of our data driven feedback control framework by solving a nonlinear optimal control problem (with non-quadratic cost functional), which can not be solved analytically. The benchmark method that we use in this example implements our data driven feedback control framework with the classic Zakai filter to calculate the conditional pdf for the controlled state  \cite{Bao_Zakai}  and the gradient decent optimization with fully calculated FBSDEs to solve the corresponding optimal control problem \cite{Gong_2017}. Both the Zakai filter and the gradient descent method are well studied with rigorous numerical analysis.
In example 3, we solve a practical application problem -- the Dubins vehicle problem, in which we control a car-like robot to reach a given target platform in the designated arrival time \cite{Marzouk_Control}.

\subsubsection*{Example 1}
Consider the following controlled process
\begin{equation}\label{Ex1:state}
dX_t = A(t)X_t dt + B U_t dt + C dW_t, \quad X_0 = x_0,
\end{equation}
where $A(t)$, $B$ and $C$ are given coefficients, $X_t$ is the multi-dimensional controlled process with the given initial position $x_0$ and $U_t$ is the multi-dimensional control process and the cost function $J$ is defined by
$$J(U) = E \Big[\frac{1}{2}\int_{0}^{T} \Big(\lan Q X_t, X_t\ran+\lan RU_t, U_t\ran \Big) dt + \frac{1}{2}\lan FX_T, X_T\ran\Big]$$
The optimal control $U_t$ for the above stochastic optimal control problem can be derived analytically as
\begin{equation}\label{Ex1:control}
\bar{U}_t = - R^{-1} B^TP(t)X_t,
\end{equation}
where $P(t)$ is the unique solution of the following Riccati equation
$$\frac{d P(t)}{dt} = - P(t) A(t) - A^T(t) P(t) + P(t) B R^{-1}  B^T P(t) - Q , \quad P(T) = F.$$
We can see from \eqref{Ex1:control} that the optimal control $\bar{U}$ depends on the state of the controlled process, which forms a feedback control problem. In our experiments, we assume that the exact state of the controlled process $X_t$ is not known and we use nonlinear indirect observations $M_t = \sin(X_t) + \eta_t$ to collect measurement data about $X_t$, where $\eta_t$ is a Gaussian random noise with covariance $\Gamma$. Therefore, the data driven feedback control has the cost functional $J^*(U)$
\begin{equation}\label{Ex1:cost}
J^*(U^M) = E \Big[\frac{1}{2}\int_{0}^{T} \Big(\lan Q X_t, X_t\ran +\lan RU^M_t, U^M_t\ran \Big) dt + \frac{1}{2}\lan F X_T, X_T\ran\Big],
\end{equation}
where $U^M$ is a data driven feedback control based on observations $M$ and we try to find $U^*$ such that $J^*(U^*) := \min_{U^M \in \mathcal{U}_{ad}} J^*(U^M)$.

In this example, we first solve the problem by choosing the controlled process as a 2D system with 1D control process. Specifically, let $A(t) = (2\sin(t), \cos(t)) \cdot I_2$ $B = (0.5, 0.5) \cdot I_2$ and $C = 0.1 I_2$ in \eqref{Ex1:state}.  For the cost functional, we let $Q = I_2$, $R= 1$ and $F = I_2$. We solve the data driven feedback control problem over the time interval $[0, 1]$ by choosing the temporal step-size $\Delta t = 0.02$, i.e. $N_T = 50$, and we let $\Gamma$ in the observational noise be $\Gamma = (0.1)^2 I_2$.  The initial state of the controlled process is chosen as $x_0 = (1, -2)^T$.
\begin{figure}[h!]
\begin{center}
\subfloat[Estimation for the control process]{\includegraphics[scale = 0.65]{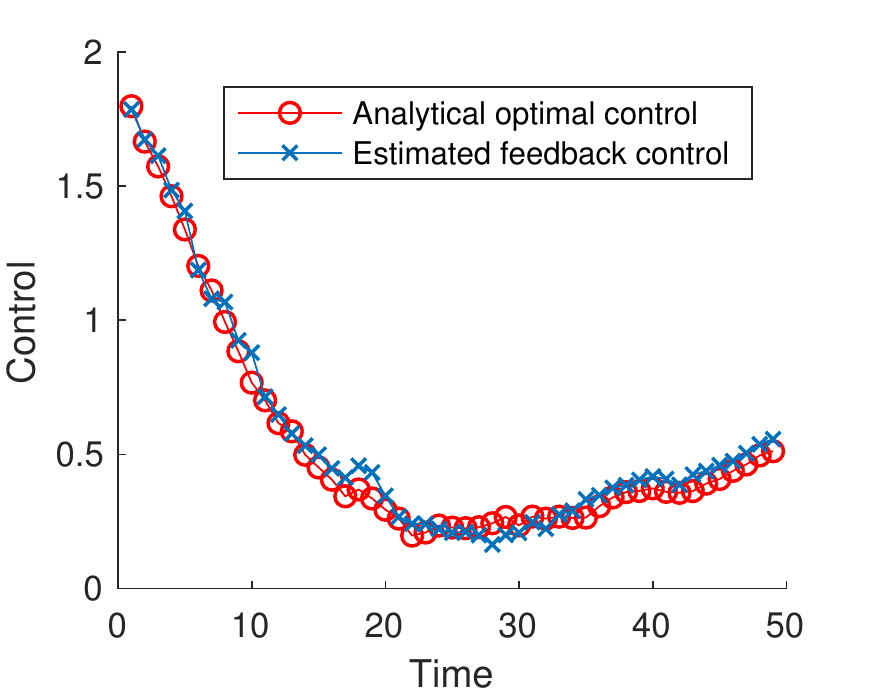} }  \quad
\subfloat[Estimation for the controlled process]{\includegraphics[scale = 0.65]{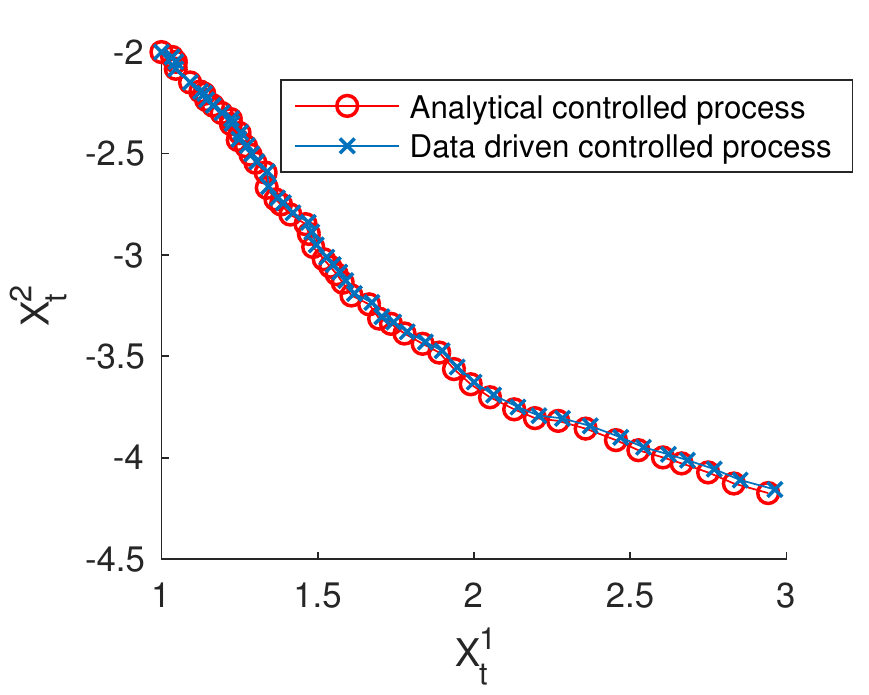} }\end{center}
\caption{Comparison between empirical representation for PDF $u$ and Shepard's approximation for PDF $u$ }\label{Ex1_2D_Performance} \vspace{-0.5em}
\end{figure}
In Figure \ref{Ex1_2D_Performance}, we demonstrate the effectiveness of our data driven feedback control algorithm by comparing with the analytical solution, where the analytical optimal control is calculated from \eqref{Ex1:control} with the \textit{exact controlled state} and we let the analytical optimal control drive its corresponding controlled system through equation \eqref{Ex1:state} with synthetically generated Brownian motion $W$ to get the analytical controlled process.
The estimated optimal control and its corresponding controlled process are calculated by our data driven feedback control \textit{Algorithm 1} based on observations $M$ with $500$ particles in the particle filter to describe the conditional pdf and $1000$ iteration steps for the SGD optimization for the optimal control, i.e. $S = 500$, $L = 1000$.  In Figure \ref{Ex1_2D_Performance} (a), we compare our estimated data driven feedback control (blue curve marked by crosses) with the analytical optimal control (red curve marked by circles); and In Figure \ref{Ex1_2D_Performance} (b), we compare the controlled process driven by our estimated feedback control  (blue curve marked by crosses) with the analytical optimal controlled process (red curve marked by circles). We can see from this figure that our algorithm accurately recovers the analytical optimal control process based on partial noisy observations for the controlled process, and the controlled process governed by our data driven feedback control well aligns with the real state of analytical controlled process.
To further examine the performance of our algorithm, we solve this 2D feedback control problem repeatedly and compare our estimated control with the analytical control to get the accumulated root mean square errors (RMSEs), i.e.
\begin{equation}\label{Ex1:RMSEs}
Err_{RMSE} : = \sqrt{\frac{1}{M_{rept}} \sum_{m=1}^{M_{rept.}}\sum_{n=1}^{N_T}\|\hat{U}^{(m)}_{t_n} - U^{(m)}_{t_n}\|_{L^2}^2},
\end{equation}
where $M_{rept.}$ is the total number of repeated experiments, $\hat{U}^{(m)}_{t_n}$ is our estimated data driven feedback control at the time instant $t_n$ for the $m$-th repeated experiment and $U^{(m)}_{t_n}$ is the analytical optimal control corresponding to the $m$-th realization of the analytical controlled process at time instant $t_n$. After repeating the above numerical experiments $100$ times with different random samples, i.e. $M_{rept.} = 100$, we obtain that $\bm{Err_{RMSE} = 0.0793}$, which is very low noting that $Err_{RMSE}$ is the accumulated RMSE combining errors of all time instants over the time interval $[0, 1]$, and this confirms the accuracy of our PF-SGD algorithm in solving data driven feedback control problem.

\begin{figure}[h!]
\begin{center}
\subfloat[$X_t^1$]{\includegraphics[scale = 0.57]{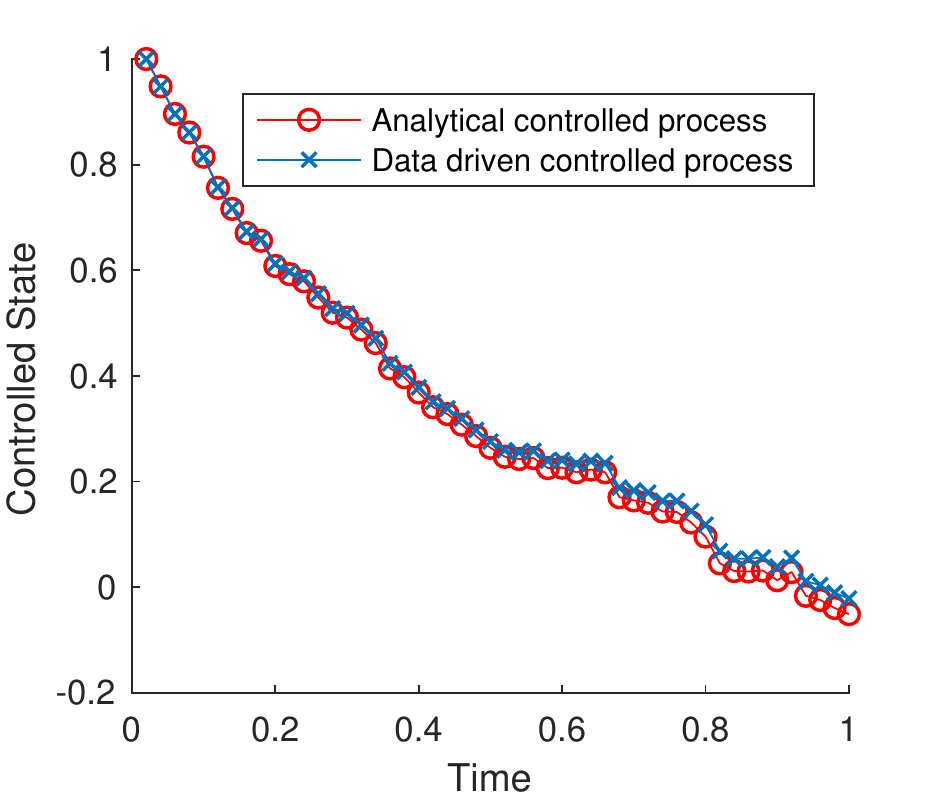} }  \quad
\subfloat[$X_t^2$]{\includegraphics[scale = 0.57]{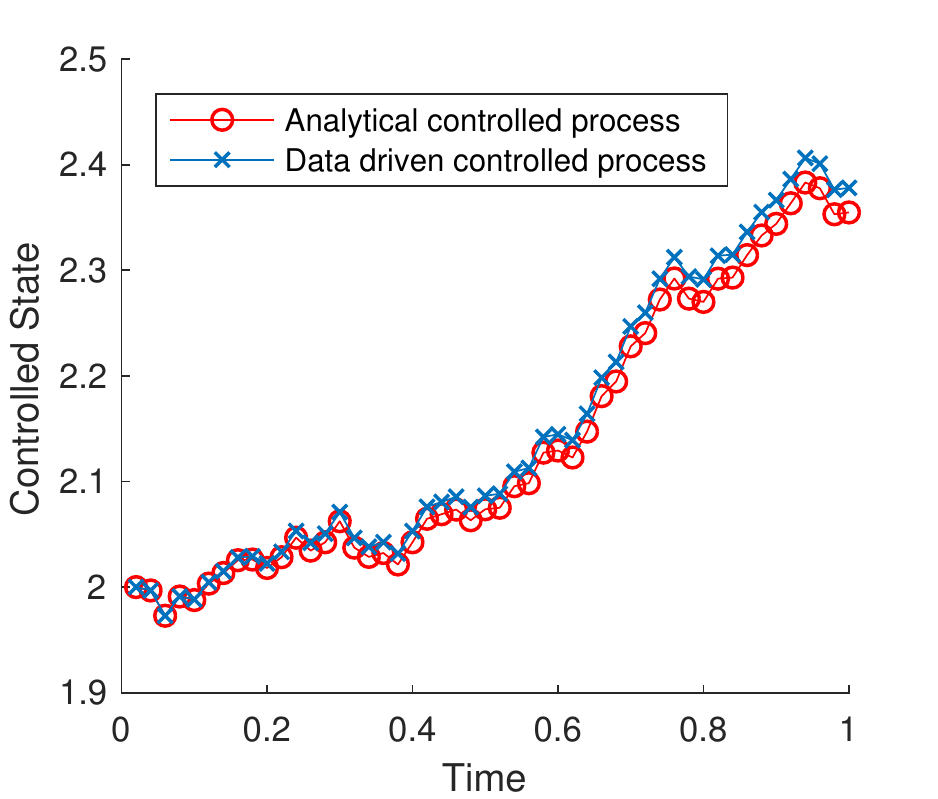} }\\
\subfloat[$X_t^3$]{\includegraphics[scale = 0.57]{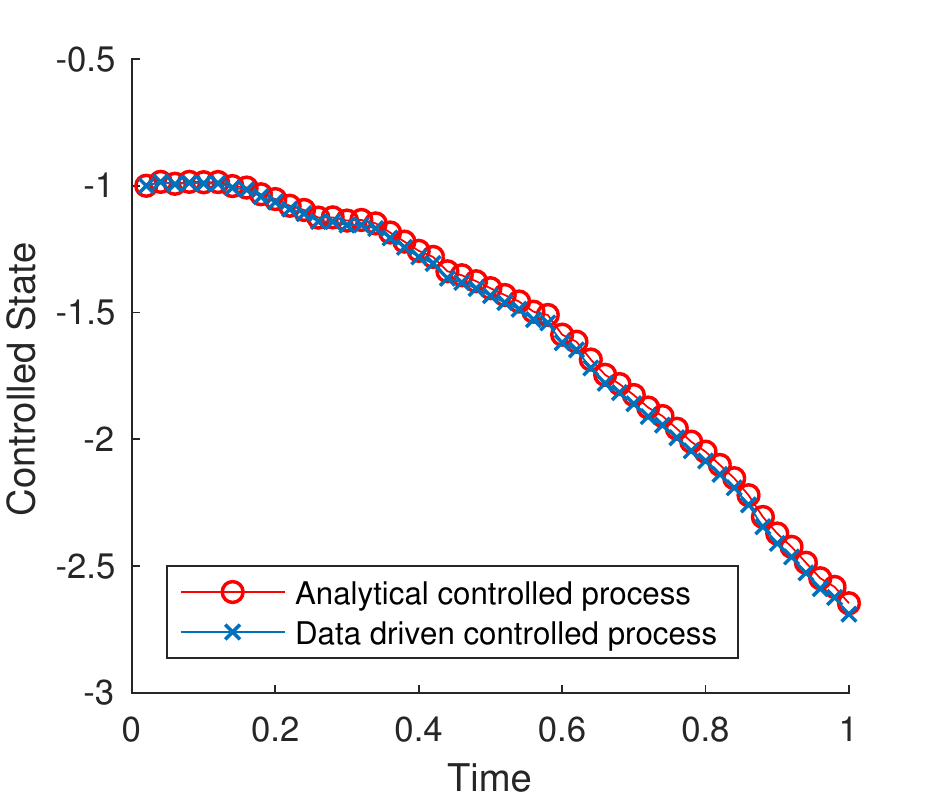} }  \quad
\subfloat[$X_t^4$]{\includegraphics[scale = 0.57]{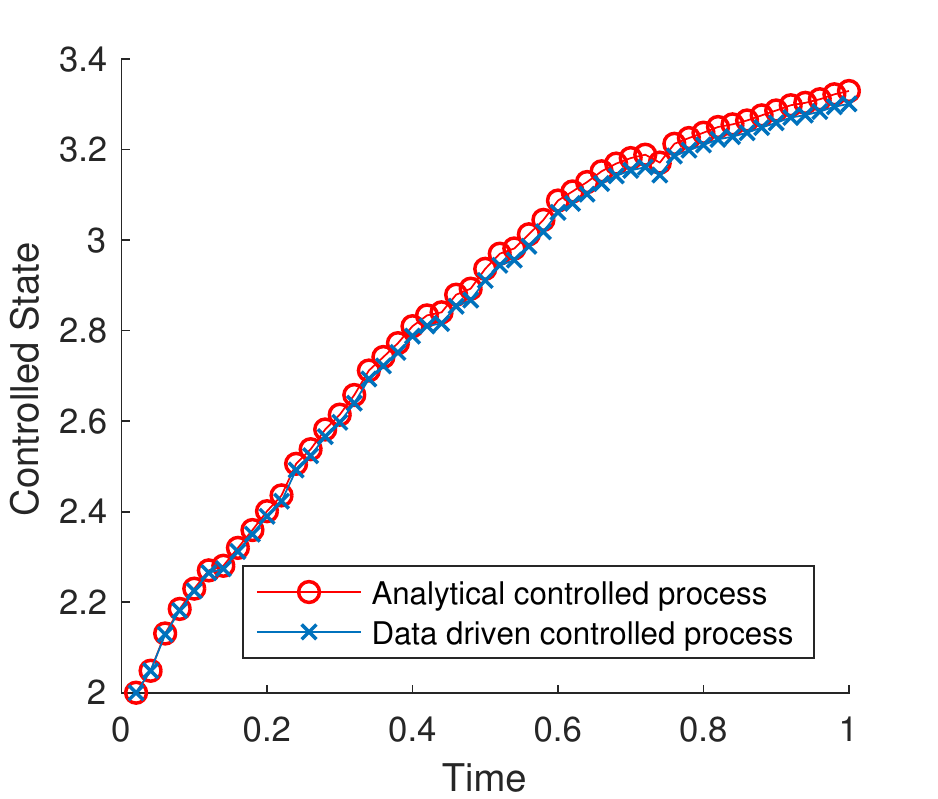} }
\end{center}
\caption{Estimation for the controlled process }\label{Ex1_4D_State} \vspace{-0.5em}
\end{figure}

One advantage of our algorithm is that it transfers the computational efforts in solving dimension-dependent stochastic differential systems to the iterations of SGD. This makes our method less sensitive to the dimension of the problem, which would address the curse of dimensionality in solving the data driven optimal control problem. In what follows, we solve a 4D version of the optimal control problem \eqref{Ex1:state} - \eqref{Ex1:cost}. Specifically, we let $A(t) = (2\sin(t), \cos(t), 1, 0.5) \cdot I_4$ $B = (0.5, 0.5, 1, 1) \cdot I_4$ and $C = 0.1 I_4$ in the controlled process.
For the cost functional, we let $Q = I_4$, $R= I_2$ and $F = I_4$. As a result, the controlled process is a 4D dynamical system and the control is a 2D process.
We also solve the data driven feedback control problem over the time interval $[0, 1]$ by choosing the temporal step-size $\Delta t = 0.02$, and we let the initial state of the controlled process be $x_0 = (1, 2, -1, 2)^T$.
\begin{figure}[h!]
\begin{center}
\includegraphics[scale = 0.7]{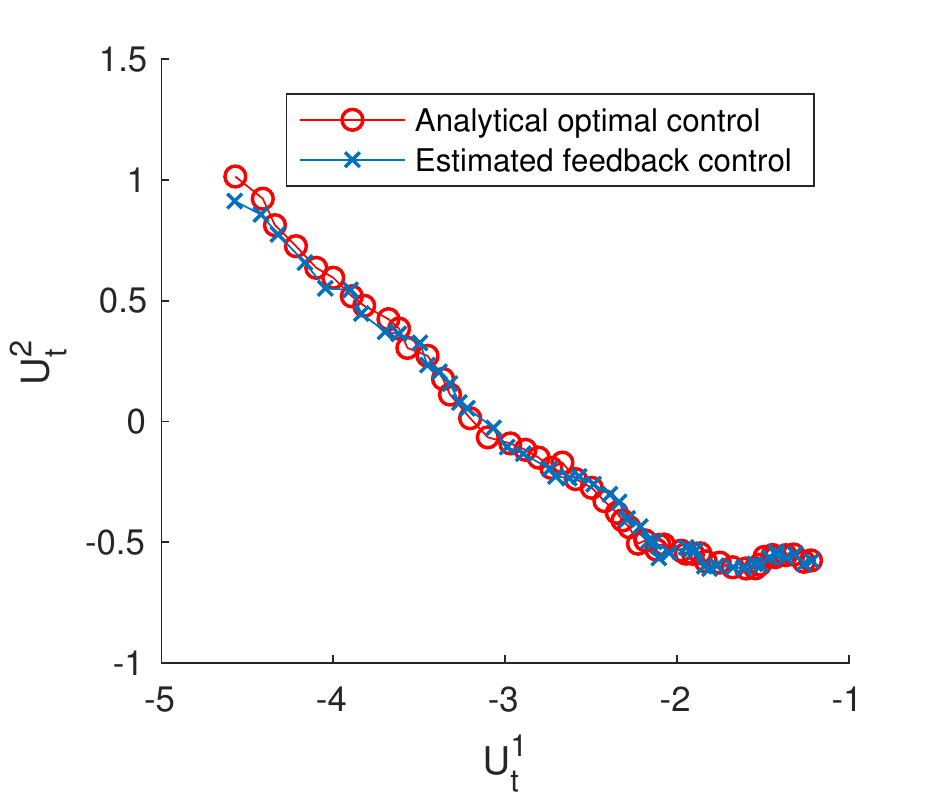}
\end{center}\vspace{-0.5em}
\caption{Estimation for the control process }\label{Ex1_4D_Control} \vspace{-.5em}
\end{figure}
Similar to the previous experiment, we compare our estimated data driven feedback control (blue curve marked by crosses) with the analytical optimal control (red curve marked by circles) in Figure \ref{Ex1_4D_Control}. To better describe a 4D distribution for the controlled state, we use $1000$ particles in the filtering stage of our algorithm and still use $1000$ iterations in the SGD optimization. In Figure \ref{Ex1_4D_State}, we compare the controlled process driven by our estimated control process (blue curve marked by crosses) with the analytical optimal controlled process (red curve marked by circles) in each dimension.
To further examine the performance of our algorithm, we repeat this 4D feedback control problem $100$ times with different random samples, i.e. $M_{rept.} = 100$. Then, we calculate the accumulated RMSEs as defined in \eqref{Ex1:RMSEs} and get $\bm{Err_{RMSE} = 0.0649}$, which combines errors on all time instants $\{t_n\}_{n=1}^{50}$.

From the above numerical experiments, we can see that our method could accurately obtain the desired feedback control in both 2D experiments and 4D experiments.
On the other hand, the average CPU time for solving the 2D problem is $9.24$ seconds and the average CPU time for solving the 4D problem is $13.51$ seconds, which is not a significant increase although both the dimension of the controlled process and the control process is doubled. Since the only difference between the 2D and 4D problems is the dimension -- while they maintain the same matrices structure, this CPU time comparison indicates that our algorithm is not very sensitive to the number of dimensions.

\subsubsection*{Example 2}
With verified effectiveness of our computational framework through numerical experiments in Example 1, in this example we focus on the demonstration of the efficiency.
Consider the following 1D nonlinear controlled process
\begin{equation}\label{Ex2:state}
\begin{aligned}
dX_t &= \arctan(X_t + u_t) dt + \sigma X_t dW_t,
\end{aligned}
\end{equation}
where $u_t$ is the control process, $\sigma X_t$ is the coefficient that determines the size of noise perturbing the controlled system and we let $\sigma = 0.05$ in this example.
The cost functional we try to minimize is
\begin{equation}\label{Ex2:cost}
J^*(u^M) = E \Big[\frac{1}{2}\int_{0}^{T}  \sin^2\big(X_t+u^M_t \big)  dt  \Big],
\end{equation}
where the observation process is defined by $M_t = X_t + \eta_t$ with a Gaussian type observational noise $\eta_t$ and we choose the standard deviation for $\eta_t$ to be $0.05$. We can see that the optimal control problem \eqref{Ex2:state}-\eqref{Ex2:cost} is a nonlinear/non-quadratic control problem, which is not explicitly solvable.

In this example, we demonstrate the efficiency of our PF-SGD algorithm by comparing with the conventional benchmark numerical implementation of the data driven feedback control framework \eqref{GD_n_s}-\eqref{J'_tn}, where the conditional pdf for the controlled process, i.e. $p(X_t|\cF^M_{t})$, is obtained by the Zakai filter through numerical solutions of the Zakai equation, and the gradient process $(J^*)_u$ is calculated by approximate solutions of the adjoint FBSDEs system. Both the Zakai equation and the FBSDEs are fully solved in the controlled state space, hence we call this conventional benchmark approach the \textit{full solution method}. We want to use this one-dimensional problem comparison to show that by avoiding to obtain full solutions of equations in each gradient descent iteration step, the PF-SGD algorithm could solve the data driven feedback control problem (C$^*$) efficiently.
Since this is a one-dimensional example, advanced high dimensional numerical techniques are not needed and we use the standard grid mesh based  linear interpolation to handle the spatial dimension approximation.
In our numerical experiments, we let $T = 1$ and implement the full solution method with a coarse mesh calculation by choosing $\Delta t = 0.1$, i.e. $N_T = 10$, with spatial partition step-size $\Delta x = 0.1$ over the pre-determined spatial interval $[3, 6]$; we also carry out a finer mesh calculation by choosing $\Delta t= 0.05$, i.e. $N_T = 20$, with spatial partition step-size $\Delta x = \frac{\sqrt{2}}{2} \cdot 0.1$ over $[3, 6]$; in the finest mesh calculation, we let $\Delta t = 0.025$, i.e. $N_T = 40$, with spatial partition step-size $\Delta x = 0.05$. In this way, for the fixed control period $[0, 1]$, we carry out more frequent control actions when choosing finer grid mesh in the full solution method.
On the other hand, the temporal step-size in our PF-SGD algorithm is $\Delta t = 0.02$, , i.e. $N_T = 50$, and we use $500$ particles to describe the conditional pdf of $X_t$ and choose the number of iteration to be $1000$ in the SGD optimization.
\begin{figure}[h!]
\begin{center}
\subfloat[Estimated control]{\includegraphics[scale = 0.65]{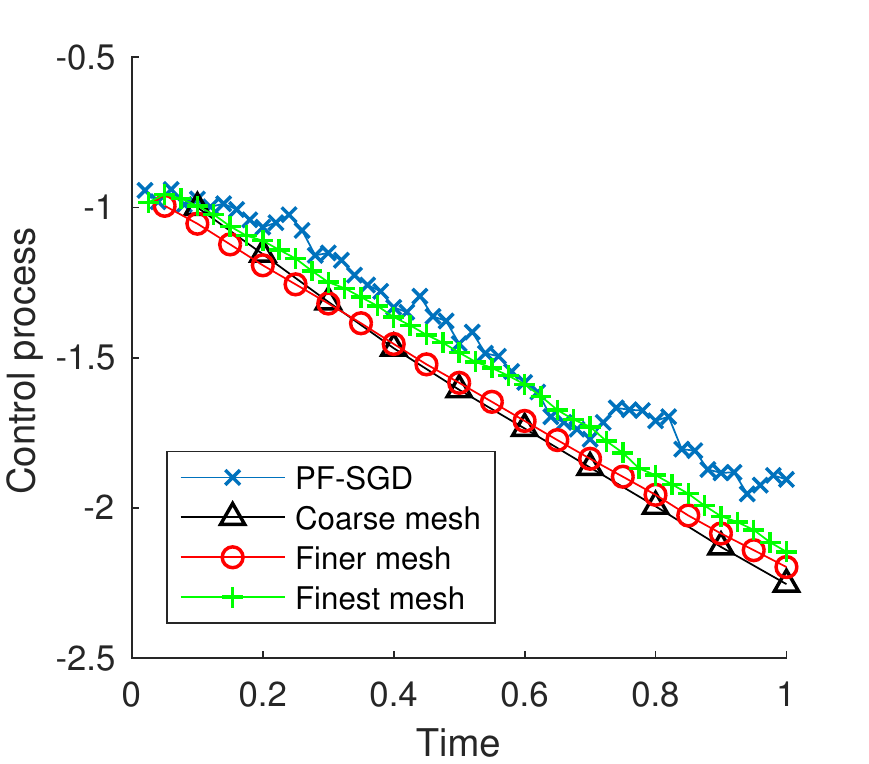} } \quad
\subfloat[Controlled process]{\includegraphics[scale = 0.65]{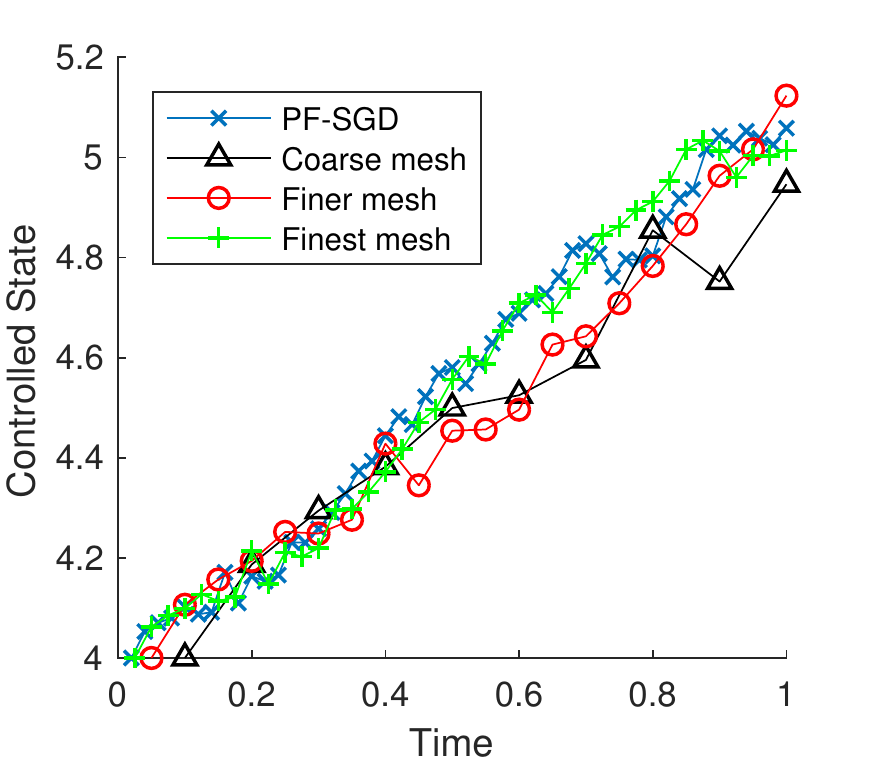}}
\end{center}
\caption{Estimation for the control process }\label{Ex2_Control-State} \vspace{-0.5em}
\end{figure}
In Figure \ref{Ex2_Control-State}, we plot the estimated controls (Figure \ref{Ex2_Control-State} (a)) and their corresponding controlled processes (Figure \ref{Ex2_Control-State} (b)). Specifically, we use the blue curve marked by ``crosses'' to represent the numerical results obtained by our PF-SGD algorithm. The black curve marked by ``triangles'' are numerical results obtained by the full solution method with coarse mesh and $10$ control actions. The red curve marked by ``circles'' are numerical results obtained by the full solution method with finer mesh and $20$ control actions. The green curve marked by ``plus signs'' are numerical results obtained by the full solution method with the finest mesh and $40$ control actions.
From this figure, we can see that the PF-SGD method could provide more detailed features for the control process which leads to more flexible actions in the controlled process; and when we refine both temporal and spatial step-sizes, the full solution method's evaluation gets closer and closer to the PF-SGD evaluation.  To indicate the performance of each algorithm, we compare the overall cost $J^*$ (as defined in \eqref{Ex2:cost}) on time interval $[0, 1]$ as well as computing time in Table \ref{Cost-sample}.
\renewcommand{\arraystretch}{1.25}
\begin{table} [h!]\small
\leftmargin=6pc \caption{} \label{Cost-sample} 
\begin{center}
\begin{tabular}{|c|c|c|c|c|}
 \hline  Method &Coarse mesh &  Finer mesh & Finest mesh & PF-SGD \\
\hline   Overall cost & $0.0481$ & $0.0318$ & $0.0076$ & $0.00095$ \\
\hline   CPU time (seconds) & $29.78$ & $220.47$ & $1560.15$ & $0.93$ \\
\hline
\end{tabular}\end{center}
\end{table}

We can see from the table that with finer temporal and spatial step-size, the full solution method gets lower cost, which can be explained by the convergence results of the Zakai filter and numerical solutions for FBSDEs, and the PF-SGD has the lowest cost. On the other hand, the full solution method spends much more CPU time when carried out with finer partition, and the PF-SGD spends very short time to reach the lowest overall cost in this experiment. The reason that the full solution method is expensive is that we have to solve both the Zakai equation and the FBSDEs system on a large domain in the controlled state space to give \textit{one update} for the estimated optimal control, and thus the computational cost for each gradient descent iteration is high.
Moreover, the corresponding computational cost would increase significantly as the dimension of the problem increases, which makes the full solution method prohibitive.

To further demonstrate the comparison results, we repeat the above experiment and calculate the average cost $\hat{J}_{t}^*(\hat{u})$ (with respect to time) defined by
$$\hat{J}_{t}^*(\hat{u}) := \frac{1}{M_{rept.}} \sum_{m=1}^{M_{rept.}}\frac{1}{2} \int_{0}^t \sin^2\big(\hat{X}_{s}^{(m)}+\hat{u}^{(m)}_{s} \big) ds,$$
where $\hat{u}^{(m)}$ is the calculated control action for the $m$-th repeated experiment and $\hat{X}^{(m)}$ is the controlled process corresponding to $\hat{u}^{(m)}$.
\begin{figure}[h!]
\begin{center}
\includegraphics[scale = 0.75]{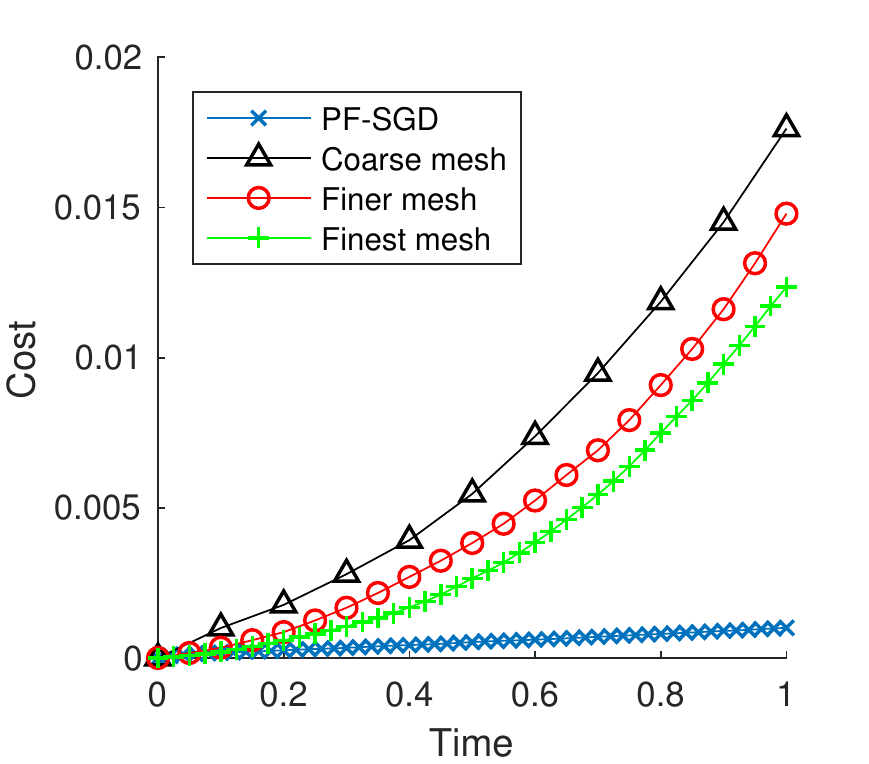}
\end{center}
\caption{Comparison of performance cost over time.}\label{Ex2_Cost}
\end{figure}
In Figure \ref{Ex2_Cost}, we present the average cost trajectories $\hat{J}_t^*(\hat{u})$ obtained from different numerical methods through $M_{rept.} = 50$ repeated experiments in solving the feedback control problem \eqref{Ex2:state}-\eqref{Ex2:cost}. The blue curve marked by ``crosses'' is the cost trajectory of our PF-SGD algorithm over time, the black curve marked by ``triangles'' is the cost trajectory of the coarse mesh implementation of the full solution method, the red curve marked by ``circles'' is the cost trajectory of finer mesh implementation of the full solution method, and the green curve marked by ``plus signs'' is the cost trajectory of the finest mesh implementation of the full solution method. We can see that when solving the problem with finer mesh and more control actions, the conventional full solution method would perform better by achieving lower cost. On the other hand, although the PF-SGD spends much less computing time, the cost trajectory of PF-SGD is much lower than full solution methods, which outperforms the conventional methods in both time and cost.

\subsubsection*{Example 3}
In this example, we solve a  Dubins vehicle maneuvering problem. The controlled process is described by the following nonlinear controlled dynamics
\begin{equation}\label{Ex3:state}
\begin{aligned}
dX_t &= \sin(\theta_t) dt + \sigma dW_t,\\
dY_t &= \cos(\theta_t) dt + \sigma dW_t,\\
d\theta_t & = u_t dt + \sigma^2 dW_t,
\end{aligned}
\end{equation}
where the pair $(X, Y)$ gives the position of a car-like robot moving in the 2D plane, $\theta$ is the steering angle that controls the moving direction of the robot, which is governed by the control action $u_t$, $\sigma$ is the noise that perturbs the motion and control actions. In our numerical experiments, we choose $\sigma = 0.2$.
The performance cost functional based on observational data that we aim to minimize is defined as
\begin{equation}\label{Ex3:cost}
J^*(u^M) = E \Big[\int_{0}^{T} \frac{1}{2} (u^M_t)^2 dt +  \delta \Big( (X_T - X_P)^2 + (Y_T - Y_P)^2 \Big)\Big],
\end{equation}
where $(X_P, Y_P)$ is a target location for the robot. The goal of this optimal control problem is to let the robot reach the target location at the given terminal time $T$ while trying to minimize the control actions, and we let $\delta = 10$ to give a strong enforcement to emphasize the importance of arriving at the target location at the terminal time. We also want to mention that the speed of the robot is an invariant constant, and the robot has to move appropriately so that it arrives at the target location at the desired time $T$.
In order to make our experiments more practical, we assume that we do not have direct observations on the robot. Instead, we use two detectors located on different observation platforms at $(6, 1)$ and $(-1, 4)$ to collect bearing angles of the target robot as indirect observations. Specifically, the observation process is defined by
$$M_t = [ \arctan(\frac{X_t - 6}{Y_t - 1}), \arctan(\frac{X_t +1}{Y_t - 4})]^T + \eta_t,$$
where the observational noise $\eta_t$ is a two dimensional Gaussian random variable with standard deviation $0.01 I_2$.

In our numerical experiments, we let $T = 1$, $\Delta t = 0.02$, i.e. $N_T = 50$, and we choose the target terminal location of the robot to be $(X_P, Y_P) = (5, 3)$. In this way, we process observational data and implement control actions $50$ times during the feedback control time period $[0, 1]$. In the PF-SGD algorithm, we use $1000$ particles to describe the conditional distribution for the controlled process and choose $1000$ iteration steps in SGD optimization.
\begin{figure}[h!] \vspace{-0.75em}
\begin{center}
\subfloat[]{\includegraphics[scale = 0.55]{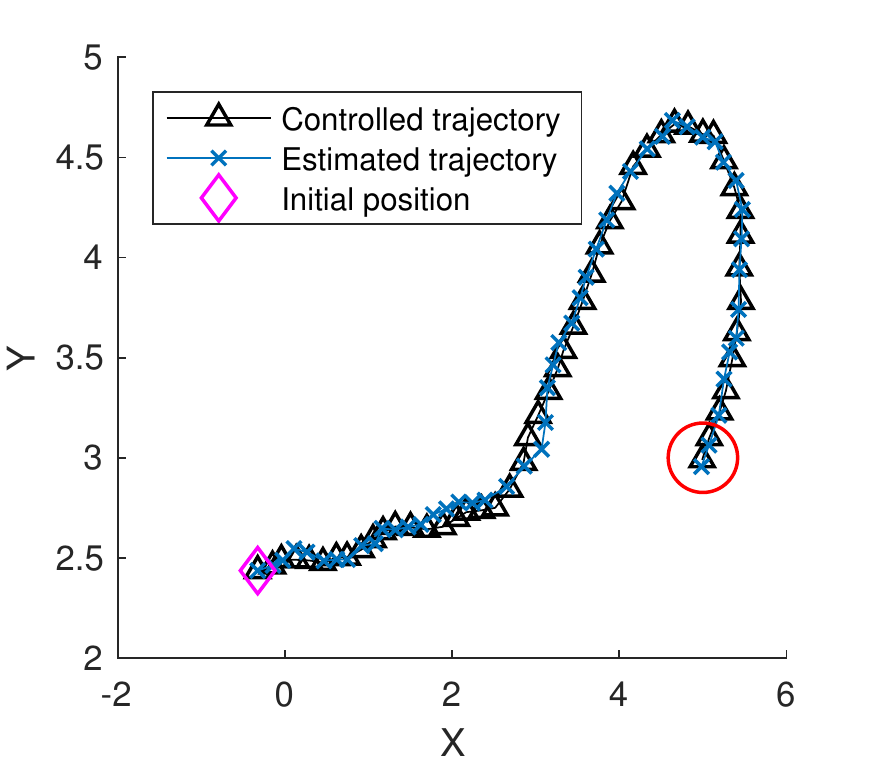} }
\subfloat[]{\includegraphics[scale = 0.55]{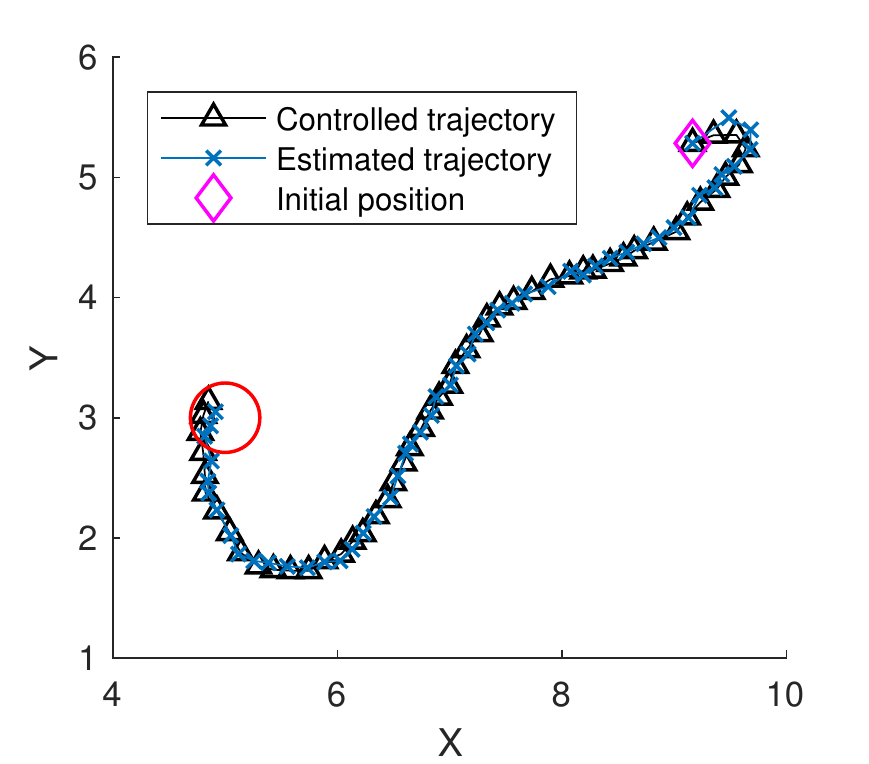}}
\subfloat[]{\includegraphics[scale = 0.55]{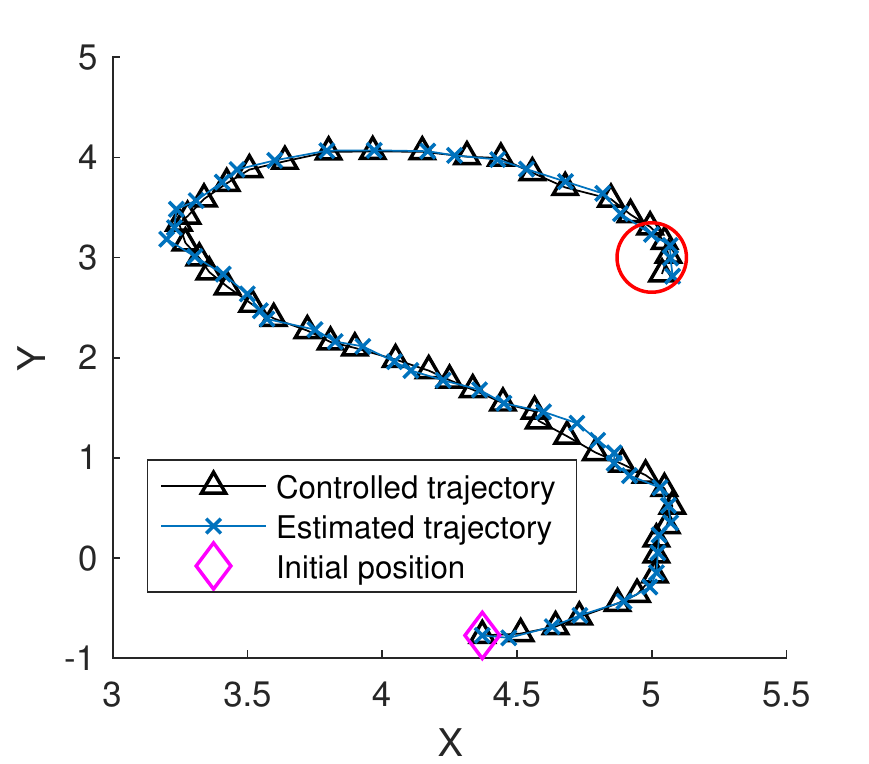}}
\end{center}
\caption{Controlled trajectories with different initial points }\label{Ex3_Initials}
\end{figure}
In Figure \ref{Ex3_Initials}, we present trajectories of the controlled car-like robot driven by observational data in three repeated experiments with three different random initial positions/steering angles. The red circle around point $(5, 3)$ is the target at the terminal time; the magenta diamonds mark different initial positions of the robot; the black curves marked by triangles give the real positions of the robot; and the blue curves marked by crosses are estimated positions of the controlled robot from the optimal filter. Due to the high effectiveness of the particle filter method in solving nonlinear filtering problems, we can see from this figure that our estimation for the position of the robot is very accurate. As the main measurement for the performance of data driven optimal control, the maneuvering robot always reaches the target location at the terminal time and when the initial position is too close to the target, the robot lingers around to wait for the \textit{designated} time $T$.

To provide the robustness of performance of our algorithm in solving this Dubins vehicle maneuvering problem, we repeat the above experiment $50$ times with initial position $(1, 1)$ for the robot, which is also perturbed by a standard Gaussian noise; and the initial heading direction is $\ds \pi/2$, which is also perturbed by a Gaussian noise with standard deviation $0.3$. The average CPU time to carry out our PF-SGD implementation of data driven feedback control over $50$ time steps is $7.47$ seconds, which shows the potential feasibility of our algorithm in real-time feedback control tasks.
\vspace{-0.3em}
\begin{figure}[h!]
\begin{center}
\includegraphics[scale = 0.6]{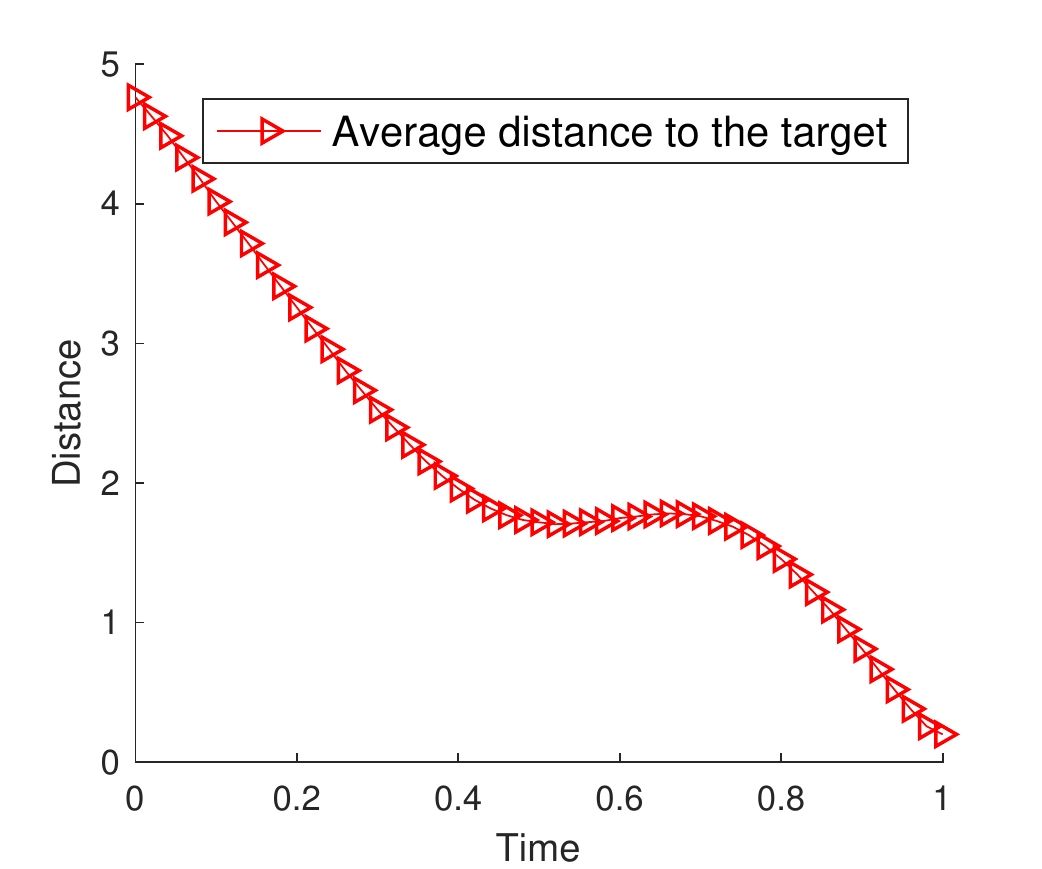}
\vspace{-0.5em}
\end{center}
\caption{Estimation for the control process }\label{Ex3_Distance} \vspace{-0.5em}
\end{figure}
In Figure \ref{Ex3_Distance}, we plot the average Euclidean distance between the robot and the target point $(5, 3)$ by the red curve marked by triangles. We can see from this figure that the robot moves towards the target platform and is ``on target'' at the final time with a very small error.

\section{Conclusions and future work}\label{Conclusion}

In this paper, we have introduced an efficient numerical algorithm to solve the data driven feedback control problem, in which the optimal control is driven by indirect observations on the (controlled) state process. The general computational framework is designed by using the particle filter method to estimate conditional expectations in the maximum principle type optimal control solver, and the optimization procedure for the control process is carried out by stochastic optimization. Numerical experiments are presented to demonstrate both the effectiveness and efficiency of our algorithm. In the future, we plan to carry out numerical analysis to derive convergence and efficiency theorems of our algorithm.  In many practical scenarios, due to the transportation of data and the computational implementation of algorithms, we also expect latency between the data generation and control actions. Therefore, we plan to incorporate the latency into the optimal control and extend the current computational framework to address delay issues in data driven feedback control problems in the future.

\newpage

\appendix

\section*{Appendix: Derivation for the gradient process}\label{Derivation for J'}

In this appendix, we give a detailed discussion on the derivation for the gradient process \eqref{J'}.
To proceed, let $U$ be convex and $(X^*,u^*)$ be any state-control pair which could be an optimal
pair. For any $u^M \in\cU_{ad}[0,T]$, we let $u^{\e,M}=u^*+\e(u^M-u^*)$.
%
%
Then the G\^ateaux derivative of $u^M\mapsto J^*(u^M)$ at $u^*$ is given by the following:
\begin{equation}\label{Variance}
\ba{ll}
 \ds \lim_{\e\to0}{J^*(u^{\e,M} )-J^*(u^*)\over\e}=\E\Big[\int_0^T\big(f^{*}_x(t) \mathcal{D}X_t+f^{*}_u(t)[u^M_t-u^{*}_t]\big)d t+h_x\cD X_T\Big],\ea
\end{equation}
where
$$\ba{ll}
\ds \cD X_t=\int_0^t\Big(b^{*}_x(s)\cD X_s+b^{*}_u(s)[u^M_s-u^*_s]\Big)ds+\int_0^t\Big(\si^{*}_x(s)\cD X_s+\si^{*}_u(s)[u^M_s-u^*_s]\Big)dW_s,\ea$$
with $\cD X_0 = 0$. For convenience of presentation, we denote $\psi^{*}(t) := \psi(t,X^*_t,u^*_t)$ and use subscript to denote partial derivative of a function. Let $(Y,Z,\z)$ be the adapted solution to the following \textit{adjoint} backward stochastic differential equation (BSDE)
$$\ds dY_s=\Big(-b^{*}_x(s)^\top Y_s-\si^{*}_x(s)^\top Z_s-g^{*}_x(s)^\top\Big)ds+Z_s dW_s+\z_sdB_s,\qq Y_T=(h^{*}_x)^\top$$
with $h^{*}_x:=h_x(X^*_T)$, where $Z$ is the martingale representation of $Y$ with respect to $W$ and $\z$ is the martingale representation of $Y$ with respect to $B$. Then
$$
\begin{aligned}
\ds h^{*}_x \mathcal{D}X_T=& \ \lan Y_T, \mathcal{D}X_T\ran-\lan Y_0, \mathcal{D}X_0\ran\\
\ds= & \  \int_0^T\Big(\lan-b^{*}_x(t)^\top Y_t-\si^{*}_x(t)^\top Z_t-f^{*}_x(t)^\top,\mathcal{D}X_t\ran+\lan Y_t,b^{*}_x(t)\mathcal{D}X_t+b^{*}_u(t)[u^M_t-u^*_t]\ran\\
& \ \ds\qq \qq +\lan Z_t,\si^{*}_x(t)\mathcal{D}X_t+\si^{*}_u(t)[u^M_t-u^*_t]\ran\Big)ds\\
& \  \ds \qq \qq +\int_0^T\Big(\lan Z_t, \mathcal{D}X_t\ran+\lan Y_t, \si^{*}_x(t)\mathcal{D}X_t+\si^{*}_u(t)[u^M_t-u^{*}_t]\Big)dW_t \\
 & \ \ds\qq\qq+\int_0^T\lan\z_t,\cD X_t+\si^{*}_u(t)[u^M_t-u_t^*]\ran dB_t\\
\ds= & \  \int_0^T\Big(-f^{*}_x(t)\mathcal{D}X_t+\lan b^{*}_u(t)^\top Y_t+\si^{*}_u(t)^\top Z_t,u^M_t-u^*_t\ran\Big)dt\\
& \  \ds\qq+\int_0^T\Big(\lan Z_t, \mathcal{D}X_t\ran+\lan Y_t,\si^{*}_x(t)\mathcal{D}X_t\ran+\si^{*}_u(t)[u^M_t-u^*_t]\Big)dW_t\\
& \  \ds\qq\qq+\int_0^T\lan\z_t,\cD X_t+\si^{*}_u(t)[u^M_t-u_t^*]\ran dB_t.
\end{aligned}
$$

Substituting the above equation into the right hand side of \eqref{Variance}, we obtain
$$
\begin{aligned}
& \lim_{\e\to0}{J^*(u^{\e, M})-J^*(u^*)\over \e}\\
=& \ \E\Big[\int_0^T\big(f_x(t)\mathcal{D}X_t+f_u(t)[u^M_t-u^*_t]\big)dt+h_x\cD X_T \Big]\\
=& \ \E\Big[\int_0^T\Big(\lan b^{*}_u(t)^\top Y_t+\si^{*}_u(t)^\top Z_t+f^{*}_u(t)^\top,u^M_t-u^*_t\ran\Big)dt\Big]\\
=& \ \E\Big[\int_0^T\lan\E\big[b^{*}_u(t)^\top Y_t+\si^{*}_u(t)^\top Z_t+f^{*}_u(t)^\top\bigm|\cF^M_t\big],u^M_t-u^*_t\ran dt\Big].
\end{aligned}
$$

Here, we have used that $u^M_t-u_t^*$ is $\cF_t^M$-measurable. Hence, in the case that $u^*$ is in the interior of $U$, one has
$$
(J^*)'_u(u^*_t)=\E\big[b^{*}_u(t)^\top Y_t+\si^{*}_u(t)^\top Z_t+f^{*}_u(t)^\top\bigm|\cF^M_t\big], \quad t \in[0, T],
$$
as required in \eqref{J'}, where $Y$ and $Z$ are solutions of the FBSDE system \eqref{FBSDE}.

\end{document}